\documentclass[12pt,a4paper,reqno]{amsart}
\usepackage[all,poly,color]{xy}
\usepackage{amsmath,amssymb, amsthm, amsfonts}
\usepackage{color}
\usepackage[pagebackref,colorlinks]{hyperref}
\usepackage{enumerate}
\usepackage{comment}

\theoremstyle{plain}
\newtheorem {lemma}{Lemma}[section]
\newtheorem {proposition}[lemma]{Proposition}
\newtheorem {theorem}[lemma]{Theorem}
\newtheorem {corollary}[lemma]{Corollary}

\theoremstyle{definition}
\newtheorem {definition}[lemma]{Definition}
\newtheorem {remark}[lemma]{Remark}
\newtheorem {example}[lemma]{Example}
\newtheorem {question}[lemma]{Question}

\newcommand{\End}{\operatorname{End}}
\newcommand{\Ann}{\operatorname{ann}}

\newcommand{\id}{\operatorname{id}}

\newcommand{\Walk}{\operatorname{Walk}}
\newcommand{\Red}{\operatorname{Red}}

\newcommand{\C}{\mathcal{C}}
\newcommand{\Es}{\mathcal{S}}
\DeclareMathOperator{\RG}{\mathbf{RG}}

\DeclareMathOperator{\Modd}{\mathbf{Mod}}
\DeclareMathOperator{\ABS}{\mathbf{ABS}}

\def\op{{\rm op}}

\begin{document}

\title[Modules for Leavitt path algebras of bi-separated graphs]{Modules for Leavitt path algebras of bi-separated graphs via representations graphs}
\author{Raimund Preusser}
\address{School of Mathematics and Statistics, Nanjing University of Information Science \& Technology, Nanjing, China} 
\email{raimund.preusser@gmx.de}

\begin{abstract}
Leavitt path algebras of bi-separated graphs have been recently introduced by R. Mohan and B. Suhas. These algebras provide a common framework for studying various generalisations of Leavitt path algebras. In this paper we obtain modules for the Leavitt path algebra $L(\dot E)$ of a finitely bi-separated graph $\dot{E}=(E,C,D)$ by introducing the notion of a representation graph for $\dot{E}$. Among these modules we find a class of simple modules. If the bi-separation on $E$ is the Cuntz-Krieger bi-separation (and hence $L(\dot{E})$ is isomorphic to the usual Leavitt path algebra $L(E)$), one recovers the celebrated Chen simple modules.
\end{abstract}

\maketitle


\section{Introduction}
Leavitt path algebras are algebras associated to directed graphs. They were introduced by G. Abrams and G. Aranda Pino in 2005 \cite{AA} and independently by P. Ara, M. Moreno and E. Pardo in 2007 \cite{AMP}. The Leavitt path algebras turned out to be a very rich and interesting class of algebras, whose studies so far have comprised almost 200 research papers. Two of these papers were written by Fields medalist E. Zelmanov \cite{zelmanov1,zelmanov2}. The field of Leavitt path algebras is a very active research area \cite{AR, arnone, BCGW, bock-sebandal-vilela, chirvasitu, crisp-macdonald, deo-chi, gupta-sen, hazrat-nam, HPS, ma-bao-li, martinez-pinedo-soler, preusser_ext, ranga_DM, vas, vas_2} with connections to functional analysis, symbolic dynamics, K-theory and noncommutative geometry. A comprehensive treatment of the subject can be found in the book~\cite{AASbook}. 

The definition of the Leavitt path algebras was inspired by the algebras $L(m,n)$ studied by W. Leavitt in the 1950's and 60's \cite{vitt56, vitt57, vitt62, vitt65}. Recall that for positive integers $m<n$ the Leavitt algebra $L(m,n)$ is universal with the property that $L(m,n)^m\cong L(m,n)^{n}$ as left $L(m,n)$-modules. 
The Leavitt path algebras embrace the algebras $L(1,n)$ but not the algebras $L(m,n),~m>1$. There have been several attempts to introduce a generalisation of the Leavitt path algebras which would cover the algebras $L(m,n),~m>1$ as well. In 2012, Ara and Goodearl introduced Leavitt path algebras of {\it separated graphs}, generalising the usual Leavitt path algebras~\cite{AG}. One recovers all of the Leavitt algebras $L(m,n)$ as corner rings of Leavitt path algebras of separated graphs. In 2013, R. Hazrat introduced Leavitt path algebras of {\it weighted graphs}, which simultaneously generalise the usual Leavitt path algebras and all Leavitt algebras $L(m,n)$ ~\cite{H-1}. In 2020, the author of this paper introduced Leavitt path algebras of {\it hypergraphs}, which cover the Leavitt path algebras of separated graphs and a subclass of the Leavitt path algebras of weighted graphs (namely the Leavitt path algebras of {\it vertex-weighted graphs}) \cite{Raimund2}. Finally in 2021, R. Mohan and B. Suhas introduced Leavitt path algebras of {\it bi-separated graphs}, which embrace all classes of algebras mentioned earlier in these paragraph \cite{mohan-suhas}.

There have been a substantial number of papers devoted to the representation theory of the usual Leavitt path algebras (of directed graphs). P. Ara and M. Brustenga proved that the category of modules for a Leavitt path algebra $L(E)$ of a finite directed graph $E$ is equivalent to a quotient of the category of modules for the path algebra $P(\overline{E})$  where $\overline{E}$ denotes the inverse graph \cite{AB}. A similar statement for graded modules over a Leavitt path algebra was established by S. Smith~\cite{smith2}. 
D. Gon\c{c}alves and D. Royer obtained modules for Leavitt path algebras by introducing the notion of an algebraic branching system \cite{GR}. X. Chen used infinite paths in $E$ to obtain simple modules for the Leavitt path algebra $L(E)$ \cite{C}. Numerous work followed, noteworthy the work of P. Ara and K. Rangaswamy producing new simple modules  and characterising those Leavitt path algebras which have countably (finitely) many distinct isomorphism classes of simple modules \cite{ARa,R-2,R-3}. G. Abrams, F. Mantese and A. Tonolo studied the projective resolutions of Chen simple modules \cite{AMT}. A recent work of P. \'Anh and T. Nam  provides another way to describe the Chen and Rangaswamy simple modules \cite{nam}. 

In \cite{GR-4}, Gon\c{c}alves and Royer obtained modules for Leavitt path algebras of separated graphs by introducing the notion of a branching system for a separated graph. In \cite{HPS}, a large class of modules for Leavitt path algebras of weighted graphs was found by introducing the notion of a representation graph for a weighted graph. For the usual Leavitt path algebras one recovers the modules via branching systems and in particular the Chen simple modules. 

In this note we describe for the first time modules for the Leavitt path algebras of bi-separated graphs. In Section 2 we recall some basic notions and results for directed graphs, as well as the definition of a bi-separated graph $\dot E=(E,C,D)$ and its Leavitt path algebra $L(\dot E)$. In Section 3, we introduce the notion of a representation graph for a given finitely bi-separated graph $\dot E$. We show that any connected component $\C$ of the category $\RG(\dot E)$ of representation graphs for $\dot E$ contains an object $(S_\C,\xi_C)$ such that the objects of $\C$ are precisely the coverings of $(S_\C,\xi_C)$ (up to isomorphism). In Section 4, we associate to any representation graph for $\dot E$ a right $L(\dot E)$-module. We prove that the module associated to $(S_\C,\xi_C)$ is simple, while the other objects of $\C$ (not isomorphic to $(S_\C,\xi_C)$) yield nonsimple modules. If the bi-separation on $E$ is the Cuntz-Krieger bi-separation and hence $L(\dot{E})$ is isomorphic to the usual Leavitt path algebra $L(E)$, one gets back the Chen simple modules. In Section 5, we prove that the module for $L(\dot E)$ defined by the universal covering $(T_\C,\zeta_C)$ of $(S_\C,\xi_C)$ is indecomposable. In Section 6, we introduce the notion of an $\dot E$-algebraic branching system, generalising Gon\c{c}alves' and Royers' notion of an $E$-algebraic branching system. We associate to any $\dot E$-algebraic branching system a right $L(\dot E)$-module and prove that the categories $\RG(\dot E)$ of representation graphs for $\dot E$ and $\ABS(\dot E)$ of $\dot E$-algebraic branching systems are equivalent.
 
\section{Preliminaries}

Throughout the paper $K$ denotes a field and $K^{\times}:=K\setminus\{0\}$. By a $K$-algebra we mean an associative (but not necessarily commutative or unital) $K$-algebra. 

\subsection{Graphs}

A {\it (directed) graph} is a quadruple $E=(E^0,E^1,s,r)$ where $E^0$ and $E^1$ are sets and $s,r:E^1\rightarrow E^0$ maps. The elements of $E^0$ are called {\it vertices} and the elements of $E^1$ {\it edges}. If $e$ is an edge, then $s(e)$ is called its {\it source} and $r(e)$ its {\it range}. In order to simplify the exposition, all graphs are assumed to contain at least one vertex and to be connected (see below), unless otherwise stated! 

Let $E$ and $F$ be graphs. $F$ is called a {\it subgraph} of $E$ if $F^0\subseteq E^0$, $F^1\subseteq E^1$, $s_F=s_E|_{F^0}$ and $r_F=r_E|_{F^0}$. A {\it graph homomorphism} $\phi: E\to F$ consists of two maps $\phi^0 : E^0\to F^0$ and $\phi^1 : E^1\to F^1$ such that $s(\phi^1(e)) = \phi^0(s(e))$ and $r(\phi^1(e)) = \phi^0(r(e))$ for any $e\in E^1$. We usually write $\phi(v)$ instead of $\phi^1(v)$ and $\phi(e)$ instead of $\phi^1(e)$ for $v\in E^0$ and $e\in E^1$. 


A {\it path} in a graph $E$ is either a finite sequence $e_1\dots e_n$ of edges such that $r(e_i)=s(e_{i+1})$ for any $1\leq i\leq n-1$, or a single vertex $v$. 
In the former case the path is called {\it nontrivial}, and in the latter case {\it trivial}. For a nontrivial path $p=e_1\dots e_n$ we set $s(p):=s(e_1)$ and $r(p):=r(e_n)$, and for a trivial path $p=v$ we set $s(p):=v$ and $r(p):=v$. The $K$-algebra $P(E)$ presented by the generating set $E^0\cup E^1$ and the relations 
\begin{enumerate}[(i)]
\item $uv=\delta_{uv}u$ for any $u,v\in E^0$ and
\smallskip
\item $s(e)e=e=er(e)$ for any $e\in E^1$
\end{enumerate}
is called the {\it path algebra} of $E$. The paths in $E$ form a linear basis for $P(E)$.

Let $E$ be a graph. The graph $\hat E=(\hat E^0, \hat E^1, \hat s, \hat r)$ where $\hat E^0=E^0$, $\hat E^1=\{e,e^*\mid e\in E^1\}$, and $\hat s(e)=s(e),~\hat r(e)=r(e),~\hat s(e^*)=r(e),~\hat r(e^*)=s(e)$ for any $e\in E^1$ is called the {\it double graph} of $E$. A path in $\hat E$ is called a {\it walk} in $E$. If $p=x_1\dots x_n$ is a walk in $E$, then $p^*:=x_n^*\dots x_1^*$ is called the {\it reverse walk} of $p$. Here we use the convention $v^*=v$ for any $v\in E^0$ and $(e^*)^*=e$ for any $e\in E^1$. We denote by $\Walk(E)$ the set of all walks in $E$. Moreover, if $u,v\in E^0$, then we denote by $_u\!\Walk(E)$ the set of all walks starting in $u$, by $\Walk_v(E)$ the set of all walks ending in $v$ and by $_u\!\Walk_v(E)$ the intersection of $_u\!\Walk(E)$ and $\Walk_v(E)$. The graph $E$ is called {\it connected} if $_u\!\Walk_v(E)\neq\emptyset$ for any $u,v\in E^0$. A graph homomorphism $\phi :E\rightarrow F$ induces a map $\Walk( E)\rightarrow \Walk( F)$, which we also denote by $\phi$. 


\subsection{Coverings of graphs}
Let $\phi: F\to E$ be a graph homomorphism. If for any $v\in F^0$ the maps $\phi^1|_{s^{-1}(v)}: s^{-1}(v)\rightarrow s^{-1}(\phi^0(v))$ and $\phi^1|_{r^{-1}(v)}: r^{-1}(v)\rightarrow r^{-1}(\phi^0(v))$ are bijective, then the pair $(F,\phi)$ is called a {\it covering} of $E$. Note that if $(F,\phi)$ is a covering of $E$, then $\phi^0$ and $\phi^1$ are surjective since we are assuming that all graphs are connected. 
Let $(F,\phi)$ and $(G,\psi)$ be coverings of $E$. A {\it morphism} $\alpha:(F,\phi)\to (G,\psi)$ is a graph homomorphism $\alpha:F\to G$ such that $\psi\circ\alpha=\phi$.

\begin{lemma}\label{lemmorcover}
Let $E$ be a graph and $\alpha:(F,\phi)\to(G,\psi)$ a morphism of coverings of $E$. Then $(F,\alpha)$ is a covering of $G$.
\end{lemma}
\begin{proof}
Let $v\in F^0$. Suppose $\alpha^1(f)=\alpha^1(f')$ for some $f,f'\in s^{-1}(v)$. Then $\psi^1(\alpha^1(f))=\psi^1(\alpha^1(f'))$ and hence $\phi^1(f)=\phi^1(f')$. Since $\phi^1|_{s^{-1}(v)}$ is injective, it follows that $f=f'$. Hence $\alpha^1|_{s^{-1}(v)}$ is injective. Let now $g\in s^{-1}(\alpha^0(v))$. Then $\psi^1(g)\in s^{-1}(\psi^0(\alpha^0(v)))=s^{-1}(\phi^0(v))$. Since $\phi^1|_{s^{-1}(v)}: s^{-1}(v)\rightarrow s^{-1}(\phi^0(v))$ is surjective, there is an $f\in s^{-1}(v)$ such that $\phi^1(f)=\psi^1(g)$ and hence $\psi^1(\alpha^1(f))=\psi^1(g)$. Clearly $s(\alpha^1(f))=\alpha^0(s(f))=\alpha^0(v)$. It follows from the injectivity of $\psi^1|_{s^{-1}(\alpha^0(v))}$ that $g=\alpha^1(f)$. Thus $\alpha^1|_{s^{-1}(v)}: s^{-1}(v)\rightarrow s^{-1}(\alpha^0(v))$ is bijective. Similarly one can show that $\alpha^1|_{r^{-1}(v)}: r^{-1}(v)\rightarrow r^{-1}(\phi^0(v))$ is bijective.
\end{proof}

A covering $(F,\phi)$ of $E$ is {\it universal} if for any covering $(G,\psi)$ of $E$ there exists a morphism $(F,\phi)\to(G,\psi)$. It is well-known that any (connected) graph $E$ has a universal covering. Below we recall the construction, cf. \cite[\S 4]{kp}. 

Let $E$ be a graph. A {\it spur} is a word of the form $ee^*$ or $e^*e$ where $e\in E^1$. A walk in $E$ is called {\it reduced} if it does not contain a spur as a subword. We denote by $\Red(E)$ (resp. ${}_u\!\Red(E)$, $\Red_v(E)$, ${}_u\!\Red_v(E)$) the set of all reduced walks in $\Walk(E)$ (resp. ${}_u\!\Walk(E)$, $\Walk_v(E)$, ${}_u\!\Walk_v(E)$). The graph $E$ is called a {\it tree}, if ${}_u\!\Red_v(E)$ has precisely one element for any $u,v\in E^0$. If $p\in \Walk(E)$, then $\underline{p}$ denotes the reduced walk in $\Red(E)$ obtained from $p$ by removing spurs until it is no longer possible (if one obtains the empty word via this procedure, then $\underline{p}=s(p)$ by definition). For any $p,q\in \Red (E)$ such that $r(p)=s(q)$ we define the reduced walk $p\ast q:=\underline{pq}$. 
If $v\in E^0$, then ${}_v\!\Red_v(E)$ is a group whose composition is given by $\ast$. This group is called the {\it fundamental group of $E$ at $v$} and is denoted by $\pi(E,v)$. The isomorphism class of $\pi(E,v)$ is independent of $v$, see \cite[\S 4]{kp}. 

Fix a vertex $v\in E^0$. We define the graph  $T = T(E, v)$ by
\begin{align*}
T^0&={}_v\!\Red(E),\\
T^1&=\{ (p,e) \in  {}_v\!\Red(E) \times  E^1 \mid r(p) = s(e) \},\\
s(p,e)&=p,\\
r(p,e)&=p\ast e.
\end{align*}
Moreover, we define a graph homomorphism $\tau=\tau(E,v): T\to E$ by $\tau^0(p)=r(p)$ and $\tau^1(p,e)=e$ for any $p\in {}_v\!\Red(E)$ and $e\in E^1$.

\begin{lemma}[{\cite[Lemmas 4.5, 4.8]{kp}}]\label{lemunicover}
Let $E$ be a graph, $v\in E^0$ and $(T,\tau)=(T(E, v),\tau(E,v))$. Then $(T,\tau)$ is a universal covering of $E$. Moreover, $T$ is a tree and the isomorphism class of $T$ is independent of the choice of the base vertex $v$.
\end{lemma}

\begin{lemma}\label{lemnospurs}
Let $(F,\phi)$ be a covering of $E$ and $p\in \Walk(F)$. Then $\phi(\underline{p})=\underline{\phi(p)}$.
\end{lemma}
\begin{proof}
If $p$ is trivial, then the assertion of the lemma is obvious. Suppose now that $p=x_1\dots x_n$ where $x_1,\dots,x_n\in \hat F^1$. Let $1\leq i\leq n-1$. It clearly suffices to show that $x_ix_{i+1}$ is a spur if and only if $\phi(x_i)\phi(x_{i+1})$ is a spur. If $x_ix_{i+1}$ is a spur then obviously $\phi(x_i)\phi(x_{i+1})$ is a spur.

Suppose now that $\phi(x_i)\phi(x_{i+1})=ee^*$ for some $1\leq i\leq n-1$ and $e\in E^1$. Then $x_i=f$ and $x_{i+1}=(f')^*$ for some $f,f'\in F^1$ such that $r(f)=r(f')$. Clearly $\phi(f)=e=\phi(f')$. It follows that $f=f'$ since $(F,\phi)$ is a covering. Thus $x_ix_{i+1}$ is a spur.

Next suppose that $\phi(x_i)\phi(x_{i+1})=e^*e$ for some $1\leq i\leq n-1$ and $e\in E^1$. Then $x_i=f^*$ and $x_{i+1}=f'$ for some $f,f'\in F^1$ such that $s(f)=s(f')$. Clearly $\phi(f)=e=\phi(f')$. It follows that $f=f'$ since $(F,\phi)$ is a covering. Thus $x_ix_{i+1}$ is a spur.
\end{proof}

We will need the lemma below in \S 3.

\begin{lemma}\label{lemunitree}
Let $E$ be a graph and $(F,\phi)$ a covering of $E$. Then $(F,\phi)$ is universal if and only if $F$ is a tree. Moreover, any universal covering of $E$ is isomorphic to $(T(E, v),\tau(E,v))$ for any $v\in E^0$.
\end{lemma}
\begin{proof}
Let $v\in E^0$. First suppose that $F$ is a tree. Since $(T,\tau):=(T(E, v),\tau(E,v))$ is universal, there is a morphism $\alpha:(T,\tau)\to (F,\phi)$. By Lemma \ref{lemmorcover}, $(T,\alpha)$ is a covering of $F$. In particular $\alpha^0$ and $\alpha^1$ are surjective. We will show that $\alpha^0$ and $\alpha^1$ are also injective. Suppose that $\alpha^0$ is not injective. Then $\alpha(u)=\alpha(w)$ for some $u\neq w\in T^0$. Since $T$ is a tree, there is a reduced walk $p\in {}_u\!\Red_w(T)$, which must be nontrivial since $u\neq w$. By the previous lemma we have $\alpha(p)\in {}_{\alpha(u)}\!\Red_{\alpha(w)}(F)$. Clearly $\alpha(p)$ is nontrivial since $p$ is nontrivial. Hence ${}_{\alpha(u)}\!\Red_{\alpha(w)}(F)$ has at least two distinct elements, $\alpha(u)$ and $\alpha(p)$, which contradicts the assumption that $F$ is a tree. Hence $\alpha^0$ is injective. Next we show that $\alpha^1$ is injective. Suppose that $\alpha(e)=\alpha(f)$ for some $e,f\in T^1$. Then $s(e)=s(f)$ since $\alpha^0$ is injective. Now it follows from the ``local" injectivity of $\alpha^1$ that $e=f$. Hence we have shown that $\alpha^0$ and $\alpha^1$ are bijective. Thus $\alpha:(T,\tau)\to (F,\phi)$ is an isomorphism and therefore $(F,\phi)$ is universal. 

Next suppose that $F$ is not a tree. Assume that $(F,\phi)$ is universal. Then there is a morphism $\alpha:(F,\phi)\to (T,\tau)$. By Lemma \ref{lemmorcover}, $(F,\alpha)$ is a covering of $T$. Since $F$ is not a tree, there are $u,w\in F^0$ such that $ {}_u\!\Red_w(F)$ contains two distinct walks, say $p$ and $q$ (note that $ {}_x\!\Red_y(F)\neq\emptyset$ for any $x,y\in F^0$ since $F$ is connected). Clearly $\alpha(p)\neq\alpha(q)$ since $(F,\alpha)$ is a covering of $T$. It follows from the previous lemma that $ {}_{\alpha(u)}\!\Red_{\alpha(w)}(T)$ contains two distinct elements, namely $\alpha(p)$ and $\alpha(q)$, which contradicts the fact that $T$ is a tree. Thus $(F,\phi)$ is not universal.

Now suppose that $(F,\phi)$ is universal. Then $F$ is a tree by the previous paragraph. It follows from the first paragraph that $(F,\phi)$ is isomorphic to $(T,\tau)$.
\end{proof}


\subsection{Bi-separated graphs and their Leavitt path algebras}
A {\it bi-separated graph} is a triple $\dot{E}=(E,C,D)$ such that  
\begin{enumerate}[(i)]
\item $E=(E^0,E^1,s,r)$ is a graph,
\smallskip
\item $C=\bigsqcup_{v\in E^0}C_v$ where $C_v$ is a partition of $s^{-1}(v)$ for any $v\in E^0$,
\smallskip
\item $D=\bigsqcup_{v\in E^0}D_v$ where $D_v$ is a partition of $r^{-1}(v)$ for any  $v\in E^0$,
\smallskip
\item $|X\cap Y|\leq 1$ for any $X\in C$ and $Y\in D$.
\end{enumerate}

In this paper all bi-separated graphs $\dot{E}=(E,C,D)$ are assumed to be {\it finitely bi-separated}, i.e. $|X|<\infty$ for any $X\in C$ and $|Y|<\infty$ for any $Y\in D$.

Let $\dot{E}=(E,C,D)$ denote a bi-separated graph. Recall from \S 2.1 that $\hat E$ denotes the double graph of $E$, and $P(E)$ the path algebra of $E$. For $X\in C$ we denote by $s(X)$ the common source of the edges in $X$. For $Y\in D$ we denote by $r(Y)$ the common range of the edges in $Y$. Moreover, for $X\in C$ and $Y\in D$ we define
\[XY=YX=\begin{cases}
e,&\text{if }X\cap Y=\{e\},\\
0,&\text{if }X\cap Y=\emptyset.
\end{cases}\]

The {\it Leavitt path algebra of $\dot E$}, denoted by $L(\dot E)$, is the quotient of $P(\hat E)$ obtained by imposing the following relations:
\begin{enumerate}[(L1)]
\item $\sum\limits_{Y\in D}(XY)(YX')^*=\delta_{XX'}s(X)\quad (X,X'\in C)$,
\medskip
\item $\sum\limits_{X\in C}(YX)^*(XY')=\delta_{YY'}r(Y)\quad(Y,Y'\in D)$.
\end{enumerate}
Note that any element of $L(\dot E)$ is a $K$-linear combination of walks in $E$.

Let $A$ be a $K$-algebra. An {\it $\dot E$-family} in $A$ is a subset $\{\sigma_v, \sigma_{e}, \sigma_{e*}\mid v\in E^0, e\in E^1\}\subseteq A$ such that\\
\vspace{-0.1cm}
\begin{enumerate}[(i)]
\item $\sigma_u\sigma_v=\delta_{uv}\sigma_u\quad(u,v\in E^0)$, 
\medskip
\item
$\sigma_{s(e)}\sigma_{e}=\sigma_{e}=\sigma_{e}\sigma_{r(e)},~\sigma_{r(e)}\sigma_{e^*}=\sigma_{e^*}=\sigma_{e^*}\sigma_{s(e)}\quad(e\in E^1)$,
\medskip
\item $\sum\limits_{Y\in D}\sigma_{XY}\sigma_{(YX')^*}=\delta_{XX'}\sigma_{s(X)}\quad (X,X'\in C)$,
\medskip
\item $\sum\limits_{X\in C}\sigma_{(YX)^*}\sigma_{XY'}=\delta_{YY'}\sigma_{r(Y)}\quad(Y,Y'\in D)$.
\end{enumerate}
Here we use the convention $\sigma_0=\sigma_{0^*}=0$. It follows from the definition of $L(\dot E)$ that there is a unique $K$-algebra homomorphism $\pi: L(\dot E)\rightarrow A$ such that $\pi(v)=\sigma_v$, $\pi(e)=\sigma_{e}$ and $\pi(e^*)=\sigma_{e^*}$ for any $v\in E^0$ and $e\in E^1$. We will refer to this as the {\it Universal Property of $L(\dot E)$}.

\section{Representation graphs for bi-separated graphs} \label{Sec2}
Throughout this section $\dot E=(E,C,D)$ denotes a fixed bi-separated graph.

\subsection{Representation graphs}
\begin{definition}\label{defwp}
A {\it representation graph} for $\dot E$ is a pair $(F,\phi)$ where $F$ is a graph and $\phi:F\rightarrow E$ a graph homomorphism such that the following holds.
\begin{enumerate}[(i)]
\item For any $w\in F^0$ and $X\in C_{\phi(w)}$ there is precisely one $f\in s^{-1}(w)$ such that $\phi(f)\in X$.

\smallskip

\item For any $w\in F^0$ and $Y\in D_{\phi(w)}$ there is precisely one $f\in r^{-1}(w)$ such that $\phi(f)\in Y$.
\end{enumerate}
\end{definition}

When visualising a representation graph $(F,\phi)$ for $\dot E$, we usually label each vertex $v\in F^0$ with $\phi(v)$ and each edge $f\in F^1$ with $\phi(f)$ (see the four examples below). Hence distinct vertices (resp. edges) in $F$ may have the same label.

\begin{example}\label{exrepgr1}
Suppose that $E$ is the graph $\xymatrix{ v\ar@(dr,ur)_{e}\ar@(dl,ul)^{f}}$, $C=C_v=\{e,f\}$ and $D=D_v=\{\{e\},\{f\}\}$. Then a representation graph $(F,\phi)$ for $\dot E$ is given by
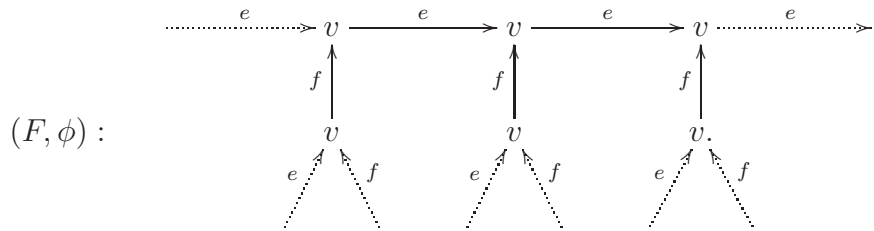
\begin{equation*}
\xymatrix@C=0.3cm{&\ar@{.>}[rrrr]^{e}&&&&v\ar[rrrr]^{e}&&&&v\ar[rrrr]^{e}&&&&v\ar@{.>}[rrrr]^{e}&&&&\\
(F,\phi):\,\,&&&&&v\ar[u]^{f}&&&&v\ar[u]^{f}&&&&v.\ar[u]^{f}&&&&\\
&&&&\ar@{.>}[ur]^{e}&&\ar@{.>}[ul]_{f}&&\ar@{.>}[ur]^{e}&&\ar@{.>}[ul]_{f}&&\ar@{.>}[ur]^{e}&&\ar@{.>}[ul]_{f}&&&}
\end{equation*}
\end{example}

\begin{example}\label{exrepgr2}
Suppose that $\dot E$ is the bi-separated graph from the previous example. Then a representation graph $(F,\phi)$ for $\dot E$ is given by
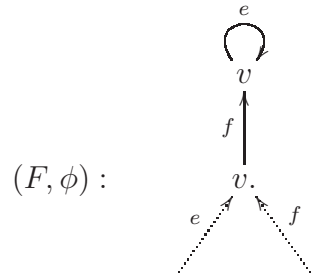
\begin{equation*}
\xymatrix@C=0.5cm{&&v\ar@(ul,ur)^{e}&\\
(F,\phi):\,\,&&v.\ar[u]^{f}&\\
&\ar@{.>}[ur]^{e}&&\ar@{.>}[ul]_{f}}
\end{equation*}
\end{example}

\begin{example}\label{exrepgr3}
Suppose that $E$ is the graph 
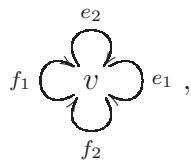
\[\xymatrix{ v\ar@(ur,dr)^{e_1}\ar@(ul,ur)^{e_2}\ar@(dl,ul)^{f_1}\ar@(dr,dl)^{f_2}},\]
$C=\{\{e_1,f_1\},\{e_2,f_2\}\}$ and $D=\{\{e_1,e_2\},\{f_1,f_2\}\}$. Then a representation graph $(F,\phi)$ for $\dot E$ is given by
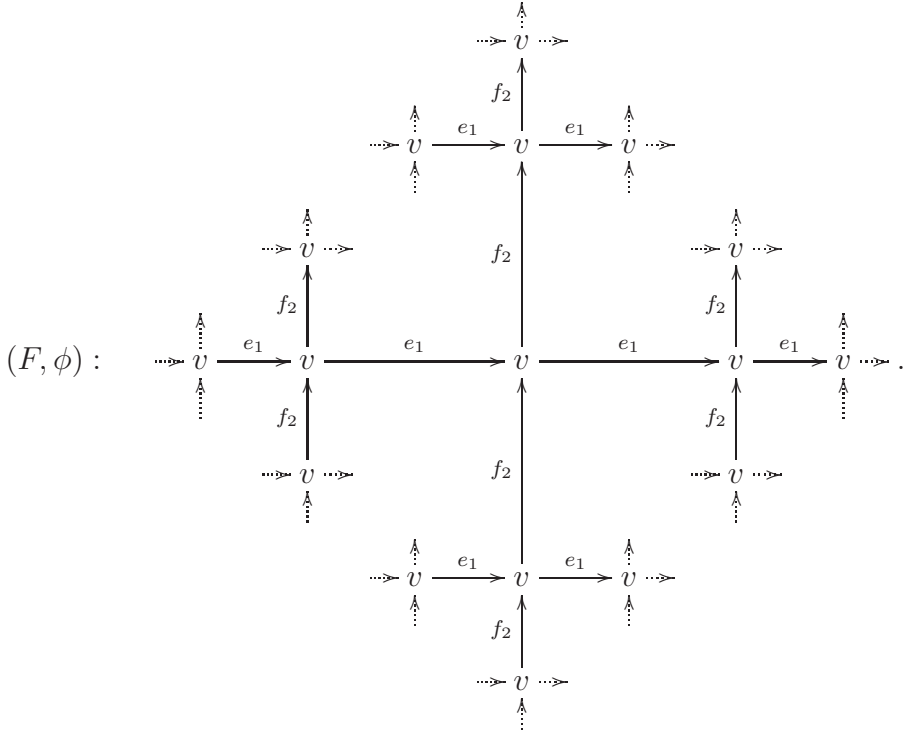
\begin{equation*}
\xymatrix@C=10pt@R=10pt{
&&&&&&&&&&&&&&&\\
&&&&&&&\ar@{..>}[r]&v\ar@{..>}[u]\ar@{..>}[r]&&&&&&&\\
&&&&&&&&&&&&&&&\\
&&&&&\ar@{..>}[r]&v\ar@{..>}[u]\ar^{e_1}[rr]&&v\ar[uu]^{f_2}\ar^{e_1}[rr]&&v\ar@{..>}[u]\ar@{..>}[r]&&&&&\\
&&&&&&\ar@{..>}[u]&&&&\ar@{..>}[u]&&&&&\\
&&&\ar@{..>}[r]&v\ar@{..>}[r]\ar@{..>}[u]&&&&&&&\ar@{..>}[r]&v\ar@{..>}[u]\ar@{..>}[r]&&&\\
&&&&&&&&&&&&&&&\\
(F,\phi):&\ar@{..>}[r]&v\ar@{..>}[u]\ar^{e_1}[rr]&&v\ar^{f_2}[uu]\ar^{e_1}[rrrr]&&&&v\ar^{f_2}[uuuu]\ar^{e_1}[rrrr]&&&&v\ar^{f_2}[uu]\ar^{e_1}[rr]&&v\ar@{..>}[u]\ar@{..>}[r]&.\\
&&\ar@{..>}[u]&&&&&&&&&&&&\ar@{..>}[u]&\\
&&&\ar@{..>}[r]&v\ar^{f_2}[uu]\ar@{..>}[r]&&&&&&&\ar@{..>}[r]&v\ar^{f_2}[uu]\ar@{..>}[r]&&&\\
&&&&\ar@{..>}[u]&&&&&&&&\ar@{..>}[u]&&&\\
&&&&&\ar@{..>}[r]&v\ar@{..>}[u]\ar^{e_1}[rr]&&v\ar^{f_2}[uuuu]\ar^{e_1}[rr]&&v\ar@{..>}[u]\ar@{..>}[r]&&&&&\\
&&&&&&\ar@{..>}[u]&&&&\ar@{..>}[u]&&&&&\\
&&&&&&&\ar@{..>}[r]&v\ar@{..>}[r]\ar^{f_2}[uu]&&&&&&&\\
&&&&&&&&\ar@{..>}[u]&&&&&&&
}
\end{equation*}
\end{example}

\begin{example}\label{exrepgr4}
Suppose that $\dot E$ is the bi-separated graph from the previous example. Then a representation graph $(F,\phi)$ for $\dot E$ is given by
\begin{equation*}
\xymatrix{(F,\phi):\,\,&v \ar@(dl,ul)^{e_1} \ar@(dr,ur)_{f_2}&.}
\end{equation*}
\end{example}

\begin{example}\label{exrepgr5}
Suppose that $E$ is the graph with one vertex $v$ and edges $e_1,e_2,f_1, g_2, h_1,h_2$. Moreover suppose that $C=\{\{e_1,f_1\},\{e_2\},\{h_1\},\{g_2,h_2\}\}$ and $D=\{\{e_1,e_2\},\{f_1\},\{g_2\},\{h_1,h_2\}\}$. Then a representation graph $(F,\phi)$ for $\dot E$ is given by
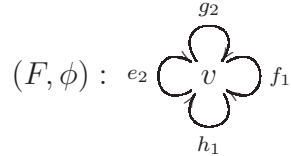
\begin{equation*}
\xymatrix{ (F,\phi):\,\,&v\ar@(ur,dr)^{f_1}\ar@(ul,ur)^{g_2}\ar@(dl,ul)^{e_2}\ar@(dr,dl)^{h_1}.}
\end{equation*}
\end{example}

The lemma below will be used quite often in the sequel.

\begin{lemma}\label{lemwell}
Let $(F,\phi)$ be a representation graph for $\dot E$. Let $q, q'\in \Walk( F)$ such that $\phi(q)=\phi(q')$. If $s(q)=s(q')$ or $r(q)=r(q')$, then $q=q'$.
\end{lemma} 
\begin{proof}
First suppose that $s(q)=s(q')=w$. If one of the walks $q$ and $q'$ is trivial, then the other must also be trivial and we have $q=w=q'$ as desired. Assume now that $q$ and $q'$ are not trivial. Then $q=x_1\dots x_n$ and $q'=y_1\dots y_n$, for some $n\geq 1$ and $x_1,\dots,x_n, y_1,\dots, y_n\in \hat F^1$. We proceed by induction on $n$.

Case $n=1$: Suppose $\phi(x_1)=\phi(y_1)=e$ for some $e\in E^1$. It follows from Definition \ref{defwp}(i) that $x_1=y_1$ and hence $q=q'$. Suppose now that $\phi(x_1)=\phi(y_1)=e^*$ for some $e\in E^1$. Then it  follows from Definition \ref{defwp}(ii) that $x_1=y_1$ and hence $q=q'$.

Case $n\to n+1$:
Suppose that $q=x_{1}\dots x_{n+1}$ and $q'=y_1\dots y_{n+1}$. By the inductive assumption we have $x_i=y_i$ for any $1\leq i\leq n$. It follows that $r(x_n)=r(y_n)$ and hence $s(x_{n+1})=s(y_{n+1})$. Now we can apply the case $n=1$ and obtain $x_{n+1}=y_{n+1}$.

Now suppose that $r(q)=r(q')$. Then $s(q^*)=s((q')^*)$. Since clearly $\phi(q^*)=\phi(q)^*=\phi(q')^*=\phi((q')^*)$, we obtain $q^*=(q')^*$. Hence $q=q'$.
\end{proof}

In \S 4 we will associate to any representation graph for $\dot E$ a module for $L(\dot E)$. The irreducible representation graphs defined below are precisely those representation graphs that yield a simple module.

\begin{definition}\label{defrepirres}
A representation graph $(F,\phi)$ for $\dot E$ is called {\it irreducible} if 
$\phi({}_u\!\Walk( F))\neq \phi({}_v\!\Walk( F))$ for any $u\neq v\in F^0$.
\end{definition}

We leave it to the reader to check that the representation graphs in Examples \ref{exrepgr1} and \ref{exrepgr3} are not irreducible, while the representation graphs in Examples \ref{exrepgr2}, \ref{exrepgr4} and \ref{exrepgr5} are irreducible.

\subsection{The category $\RG(\dot E)$}\label{catrepa} 

We denote by $\RG(\dot E)$ the category whose objects are the representation graphs for $\dot E$. A morphism $\alpha:(F,\phi)\to (G,\psi)$ in $\RG(\dot E)$ is a graph homomorphism $\alpha:F\to G$ such that $\psi\circ\alpha=\phi$. We will see that if $\alpha:(F,\phi)\to (G,\psi)$ is a morphism in $\RG(\dot E)$, then $(F,\alpha)$ is a covering of $G$. 


\begin{lemma}\label{lempath}
Let $(F,\phi)$ and $(G,\psi)$ be objects in $\RG(\dot E)$. Let $u\in F^0$ and $v\in G^0$. If $\phi({}_u\!\Walk( F))\subseteq \psi({}_v\!\Walk( G))$, then $\phi({}_u\!\Walk( F))=\psi({}_v\!\Walk( G))$.
\end{lemma}
\begin{proof}
Suppose that $\phi({}_u\!\Walk( F))\subseteq \psi({}_v\!\Walk( G))$. It follows that $\phi(u)=\psi(v)$ since $u\in {}_u\!\Walk( F)$. We have to show that $\psi({}_v\!\Walk( G))\subseteq \phi({}_u\!\Walk( F))$. Let $p\in{}_v\!\Walk( G)$. If $p=v$, then $\psi(p)=\psi(v)=\phi(u)\in \phi({}_u\!\Walk( F))$. Suppose now that $p$ is nontrivial. Then $p=y_1\dots y_n$ for some $y_1,\dots, y_n\in \hat G^1$ where $n\geq 1$. We proceed by induction on $n$.

Case $n=1$: Suppose that $p=g$ for some $g\in s_G^{-1}(v)$. Then $\psi(g)\in s^{-1}(\psi(v))$. Let $X\in C_{\psi(v)}$ such that $\psi(g)\in X$. Since $(F,\phi)$ satisfies Condition (i) in Definition \ref{defwp}, there is precisely one $f\in s_F^{-1}(u)$ such that $\phi(f)\in X$. Since $\phi({}_u\!\Walk( F))\subseteq \psi({}_v\!\Walk( G))$, there is a $g'\in s_G^{-1}(v)$ such that $\psi(g')=\phi(f)\in X$. It follows that $g=g'$ since $(G,\psi)$ satisfies Condition (i) in Definition \ref{defwp}. Thus $\psi(g)=\phi(f)\in \phi({}_u\!\Walk( F))$.

Suppose now that $p=g^*$ for some $g\in r_G^{-1}(v)$. Then $\psi(g)\in r^{-1}(\psi(v))$. Let $Y\in D_{\psi(v)}$ such that $\psi(g)\in Y$. Since $(F,\phi)$ satisfies Condition (ii) in Definition \ref{defwp}, there is precisely one $f\in r_F^{-1}(u)$ such that $\phi(f)\in Y$. Since $\phi({}_u\!\Walk( F))\subseteq \psi({}_v\!\Walk( G))$, there is a $g'\in r_G^{-1}(v)$ such that $\psi((g')^*)=\phi(f^*)$, i.e. $\psi(g')=\phi(f)\in Y$. It follows that $g=g'$ since $(G,\psi)$ satisfies Condition (ii) in Definition \ref{defwp}. Thus $\psi(g^*)=\phi(f^*)\in \phi({}_u\!\Walk( F))$.

Case $n\rightarrow n+1$:  Suppose $p=y_1\dots y_ny_{n+1}$. By the induction assumption we know that $\psi(y_1\dots y_n)\in \phi({}_u\!\Walk( F))$. Hence $\psi(y_1\dots y_n)=\phi(x_1\dots x_n)$ for some walk $x_1\dots x_n\in {}_u\!\Walk( F)$. Set $u':=r_{\hat F}(x_n)$ and $v':=r_{\hat G}(y_n)$. Clearly $\phi(u')=\psi(v')$ since $\phi$ and $\psi$ are graph homomorphisms.

Suppose that $y_{n+1}=g$ for some $g\in s_G^{-1}(v')$. Then $\psi(g)\in s^{-1}(\psi(v'))$. Let $X\in C_{\psi(v')}$ such that $\psi(g)\in X$. Since $(F,\phi)$ satisfies Condition (i) in Definition \ref{defwp}, there is precisely one $f\in s_F^{-1}(u')$ such that $\phi(f)\in X$. Since $\phi({}_u\!\Walk( F))\subseteq \psi({}_v\!\Walk( G))$, there is a $g'\in s_G^{-1}(v')$ such that $\phi(x_1\dots x_nf)=\psi(y_1\dots y_ng')$, which implies $\psi(g')=\phi(f)\in X$. It follows that $g=g'$ since $(G,\psi)$ satisfies Condition (i) in Definition \ref{defwp}. Thus $\psi(g)=\phi(f)$ and hence $\psi(y_1\dots y_ng)=\phi(x_1\dots x_nf)\in \phi({}_u\!\Walk( F))$.

Suppose now that $y_{n+1}=g^*$ for some $g\in r_G^{-1}(v')$. Then $\psi(g)\in r^{-1}(\psi(v'))$. Let $Y\in D_{\psi(v')}$ such that $\psi(g)\in Y$. Since $(F,\phi)$ satisfies Condition (ii) in Definition \ref{defwp}, there is precisely one $f\in r_F^{-1}(u')$ such that $\phi(f)\in Y$. Since $\phi({}_u\!\Walk( F))\subseteq \psi({}_v\!\Walk( G))$, there is a $g'\in r_G^{-1}(v')$ such that $\phi(x_1\dots x_nf^*)=\psi(y_1\dots y_n(g')^*)$, which implies $\psi(g')=\phi(f)\in Y$. It follows that $g=g'$ since $(G,\psi)$ satisfies Condition (ii) in Definition \ref{defwp}. Thus $\psi(g^*)=\phi(f^*)$ and hence $\psi(y_1\dots y_ng^*)=\phi(x_1\dots x_nf^*)\in \phi({}_u\!\Walk( F))$.
\end{proof}
 
\begin{lemma}\label{propcruc}
Let $\alpha:(F,\phi)\to (G,\psi)$ be a morphism in $\RG(\dot E)$. If $u\in F^0$, then $\phi({}_u\!\Walk( F))=\psi({}_{\alpha(u)}\!\Walk( G))$.
\end{lemma}
\begin{proof}
Clearly $\phi({}_u\!\Walk( F))\subseteq\psi({}_{\alpha(u)}\!\Walk( G))$ since $\phi=\psi\circ \alpha$	. It follows from Lemma \ref{lempath} that $\phi({}_u\!\Walk( F))=\psi({}_{\alpha(u)}\!\Walk( G))$.
\end{proof}

\begin{proposition}\label{propsurj}
Let $\alpha:(F,\phi)\to (G,\psi)$ be a morphism in $\RG(\dot E)$. Then $(F,\alpha)$ is a covering of $G$. 
\end{proposition}
\begin{proof}


Let $u\in F^0$. First we show that $\alpha^1|_{s^{-1}(u)}: s^{-1}(u)\rightarrow s^{-1}(\alpha(u))$ is injective. Let $f\neq f'\in s^{-1}(u)$. Let $X,X'\in C_{\phi(u)}$ such that $\phi(f)\in X$ and $\phi(f')\in X'$. Then $X\neq X'$ by the definition of a representation graph. It follows that $\alpha(f) \not = \alpha(f')$ since $\psi(\alpha(f))=\phi(f)\in X$ but $\psi(\alpha(f'))=\phi(f')\in X'$.

Next we show that $\alpha^1|_{s^{-1}(u)}: s^{-1}(u)\rightarrow s^{-1}(\alpha(u))$ is surjective. Let $g\in s^{-1}(\alpha(u))$. It follows from Lemma \ref{propcruc} that there is an $f\in s^{-1}(u)$ such that $\phi(f)=\psi(g)$. Clearly $\alpha(f)\in s^{-1}(\alpha(u))$ and $\psi(\alpha(f))=\phi(f)=\psi(g)$. It follows from the definition of a representation graph that $\alpha(f)=g$.

A similar argument shows that $\alpha^1|_{r^{-1}(u)}: r^{-1}(u)\rightarrow r^{-1}(\alpha(u))$ is bijective. 
\end{proof}


\begin{proposition}\label{propsurj2}
If $(G,\psi)$ is an object in $\RG(\dot E)$ and $(F,\alpha)$ is a covering of $G$, then $(F,\psi\circ\alpha)$ is an object in $\RG(\dot E)$ (and hence $\alpha:(F,\psi\circ\alpha)\to (G,\psi)$ is a morphism in $\RG(\dot E)$). 
\end{proposition}
\begin{proof}
We have to show that conditions (i) and (ii) in Definition \ref{defwp} are satisfied.
\begin{enumerate}[(i)]
\item Let $v\in F^0$ and $X\in C_{\psi(\alpha(v))}$. Since $(G,\psi)$ is a representation graph, there is precisely one $g\in s^{-1}(\alpha(v))$ such that $\psi(g)\in X$. It follows that there is precisely one $f\in s^{-1}(v)$ such that $\psi(\alpha(f))\in X$ since $(F,\alpha)$ is a covering.

\smallskip

\item Let $v\in F^0$ and $Y\in D_{\psi(\alpha(v))}$. Since $(G,\psi)$ is a representation graph, there is precisely one $g\in r^{-1}(\alpha(v))$ such that $\psi(g)\in Y$. It follows that there is precisely one $f\in r^{-1}(v)$ such that $\psi(\alpha(f))\in Y$ since $(F,\alpha)$ is a covering.
\end{enumerate}
\end{proof}

\subsection{Quotients of representation graphs}
For any object $(F,\phi)$ in $\RG(\dot E)$ we define an equivalence relation $\sim$ on $F^0$ by $u\sim v$ if $\phi({}_u\!\Walk( F))=\phi({}_v\!\Walk( F))$. Recall that if $\sim$ and $\approx$ are equivalence relations on a set $X$, then one writes $\approx~\leq~ \sim $ and calls $\approx$ {\it finer} than $\sim$, if $x\approx y$ implies $x\sim y$ for any $x,y\in X$. 

\begin{definition}\label{defadm}
Let $(F,\phi)$ be an object in $\RG(\dot E)$. An equivalence relation $\approx$ on $F^0$ is called {\it admissible} if the following hold:
\begin{enumerate}[(i)]
\item $\approx~\leq~\sim$.
\item If $p,q\in{} \Walk( F)$ such that $\phi(p)=\phi(q)$ and $s(p)\approx s(q)$, then $r(p)\approx r(q)$.
\end{enumerate}
\end{definition}

The admissible equivalence relations on $F^0$ with partial order $\leq$ form a bounded lattice whose maximal element is $\sim$ and whose minimal element is the equality relation $=$.

Let $(F,\phi)$ be an object in $\RG(\dot E)$ and $\approx$ an admissible equivalence relation on $F^0$. If $f, g\in F^1$ we write $f\approx g$ if $s(f)\approx s(g)$ and $\phi(f)=\phi(g)$. This defines an equivalence relation on $F^1$. Define an object $(F_\approx,\phi_\approx)$ in $\RG(\dot E)$ by 
\begin{align*}
F_\approx^0&=F^0/\approx,\\
F_\approx^1&=F^1/\approx,\\
s([f])&=[s(f)],\\
r([f])&=[r(f)],\\
\phi_\approx^0([v])&=\phi^0(v),\\
\phi_\approx^1([f])&=\phi^1(f).
\end{align*}
We call $(F_\approx,\phi_\approx)$ a {\it quotient} of $(F,\phi)$. 

\begin{lemma}\label{lemcompare}
Let $(F,\phi)$ be an object in $\RG(\dot E)$. Let $\approx~\leq ~\approx'$ be admissible equivalence relations on $F^0$. Then there is a morphism $(F_\approx,\phi_\approx)\to (F_{\approx'},\phi_{\approx'})$.
\end{lemma}
\begin{proof}
Define a graph homomorphism $\alpha:F_\approx\to F_{\approx'}$ by $\alpha^0([v]_\approx)=[v]_{\approx'}$ and $\alpha^1([f]_\approx)=[f]_{\approx'}$ for any $v\in F^0$ and $f\in F^1$. Since $\approx~\leq ~\approx'$, $\alpha$ is well-defined. Clearly $\phi_{\approx'}\circ\alpha=\phi_{\approx}$ and therefore $\alpha:(F_\approx,\phi_\approx)\to (F_{\approx'},\phi_{\approx'})$ is a morphism in $\RG(\dot E)$. 
\end{proof}

\begin{proposition}\label{thmadm}
Let $(F,\phi)$ and $(G,\psi)$ be objects in $\RG(\dot E)$. Then there is a morphism $\alpha:(F,\phi)\to (G,\psi)$ if and only if $(G,\psi)$ is isomorphic to a quotient of $(F,\phi)$.
\end{proposition}
\begin{proof}
$(\Rightarrow)$ Suppose there is a morphism $\alpha:(F,\phi)\to (G,\psi)$. If $u,v\in F^0$, we write $u\approx v$ if and only if $\alpha(u)=\alpha(v)$. Clearly $\approx$ defines an equivalence relation on $F^0$. Below we check that $\approx$ is admissible. 
\begin{enumerate}[(i)]
\item Suppose $u\approx v$. Then $\phi({}_u\!\Walk( F))=\psi({}_{\alpha(u)}\!\Walk( G))=\psi({}_{\alpha(v)}\!\Walk( G))=\phi({}_v\!\Walk( F))$ by Lemma \ref{propcruc}. Hence $u\sim v$.
\item Suppose $p\in{} _u\!\Walk_x( F)$, $q\in{} _v\!\Walk_y( F)$, $\phi(p)=\phi(q)$ and $u\approx v$. Clearly $\alpha(p)\in{} _{\alpha(u)}\!\Walk_{\alpha(x)}( G)$ and $\alpha(q)\in{} _{\alpha(v)}\!\Walk_{\alpha(y)}( G)$. Moreover, $\psi(\alpha(p))=\phi(p)=\phi(q)=\psi(\alpha(q))$. Since $\alpha(u)=\alpha(v)$, it follows from Lemma \ref{lemwell} that $\alpha(p)=\alpha(q)$. Hence $\alpha(x)=r(\alpha(p))=r(\alpha(q))=\alpha(y)$ and therefore $x\approx y$.
\end{enumerate}

Note that by Lemma \ref{lemwell} we have $f\approx g$ if and only if $\alpha^1(f)=\alpha^1(g)$, for any $f,g\in F^1$. Define a graph homomorphism $\beta:F_\approx\to G$ by $\beta^0([v])=\alpha^0(v)$ and $\beta^1([f])=\alpha^1(f)$. Clearly $\psi\circ \beta=\phi_{\approx}$ and therefore $\beta:(F_\approx,\phi_\approx)\to (G,\psi)$ is a morphism in $\RG(\dot E)$. It follows from the definition of $\beta$ that $\beta^0$ and $\beta^1$ are injective. In view of Proposition \ref{propsurj}, $\beta^0$ and $\beta^1$ are also surjective and hence $\beta$ is an isomorphism.

$(\Leftarrow)$ Suppose now that $(G,\psi)\cong (F_\approx,\phi_\approx)$ for some admissible equivalence relation $\approx$ on $F^0$. In order to show that there is a morphism $\alpha:(F,\phi)\to (G,\psi)$ it suffices to show that there is a morphism $\beta:(F,\phi)\to (F_\approx,\phi_\approx)$. But this is obvious (define $\beta^0(v)=[v]$ and $\beta^1(f)=[f]$).
\end{proof}

\subsection{The connected components of the category $\RG(\dot E)$} \label{subsecnn}
Recall that any category $\mathbf{C}$ can be written as a disjoint union (or coproduct) of a collection of connected categories, which are called the {\it connected components} of $\mathbf{C}$. Each connected component is a full subcategory of $\mathbf{C}$.

\begin{lemma}\label{lemrepell}
Let $(F,\phi)$ and $(G,\psi)$ be objects of $\RG(\dot E)$ and suppose there is a morphism $(F,\phi)\to(G,\psi)$ or a morphism $(G,\psi)\to(F,\phi)$. Let $(T,\tau)$ be a universal covering of $F$. Then there is a morphism $\eta:(T,\phi\circ \tau)\to(G,\psi)$ in $\RG(\dot E)$ such that $(T,\eta)$ is a universal covering of $G$.
\end{lemma}
\begin{proof}
Note that $T$ is a tree by Lemma \ref{lemunitree}, and that $(T,\phi\circ\tau)$ is an object in $\RG(\dot E)$ by Proposition \ref{propsurj2}. First suppose that there is a morphism $\alpha:(F,\phi)\rightarrow(G,\psi)$. Set $\eta:=\alpha\circ \tau$. Since the diagram
\[\xymatrix@C=15pt@R=20pt{&T\ar[dl]_{\tau}\ar[dr]^{\eta}&\\F\ar[dr]_{\phi}\ar[rr]^{\alpha}&&G\ar[dl]^{\psi}\\&E&}\]
commutes, $\eta:(T,\phi\circ \tau)\to(G,\psi)$ is a morphism in $\RG(\dot E)$. It follows from Proposition \ref{propsurj}, that $(T,\eta)$ is a covering of $G$. By Lemma \ref{lemunitree} this covering is universal. 
\\
Suppose now that there is a morphism $\alpha:(G,\psi)\rightarrow (F,\phi)$. Let $(T',\tau')$ be a universal covering of $G$. Then $(T',\alpha\circ \tau')$ is a universal covering of $F$ by the previous paragraph. It follows from Lemma \ref{lemunitree} that $(T, \tau)\cong (T', \alpha\circ \tau')$, i.e. there is a graph isomorphism $\gamma:T\to T'$ making the diagram
\[\xymatrix@C=15pt@R=20pt{T\ar[d]_{\tau}\ar[rr]^{\gamma}&&T'\ar[d]^
{\tau'}\\F\ar[dr]_{\phi}&&G\ar[dl]^{\psi}\ar[ll]_{\alpha}\\&E&}\]
commute. Set $\eta:=\tau'\circ\gamma$. Then clearly $\eta:(T,\phi\circ \tau)\to(G,\psi)$ is a morphism in $\RG(\dot E)$. It follows from Proposition \ref{propsurj}, that $(T,\eta)$ is a covering of $G$. By Lemma \ref{lemunitree} this covering is universal. 
\end{proof}

Let $\C$ be a connected component of $\RG(\dot E)$ and choose an object $(F,\phi)=(F_\C,\phi_\C)$ of $\C$. By \S 2.2 we can choose a universal covering $(T,\tau)$ of $F$. By Proposition \ref{propsurj2}, $(T,\phi\circ\tau)$ is an object of $\RG(\dot E)$. Obviously $\tau: (T,\phi\circ\tau)\to (F,\phi)$ is a morphism in $\RG(\dot E)$ and hence $(T,\phi\circ\tau)$ belongs to the same connected component $\C$. We set 
\[(T_\C,\zeta_\C):=(T,\phi\circ\tau)~\text{ and }~(S_\C,\xi_\C):=((T_\C)_\sim,(\zeta_\C)_\sim).\]
For the definition of the quotient $((T_\C)_\sim,(\zeta_\C)_\sim)$ see \S 3.3.

We call an object $X$ in a category $\mathbf{C}$ {\it repelling} (resp. {\it attracting}) if for any object $Y$ in $\mathbf{C}$ there is a morphism $X\rightarrow Y$ (resp. $Y\rightarrow X$). 

\begin{theorem}\label{thmm1}
Let $\C$ be a connected component of $\RG(\dot E)$. Then $(T_\C,\zeta_\C)$ is a repelling object of $\C$, and consequently the objects of $\C$ are up to isomorphism precisely the quotients of $(T_\C,\zeta_\C)$.
\end{theorem}
\begin{proof}
Let $(F,\phi)=(F_\C,\phi_\C)$ and $(G,\psi)$ be any object in $\C$. Then there is a sequence of objects \[(F,\phi)=(F_0,\phi_0),(F_1,\phi_1),\dots,(F_{n-1},\phi_{n-1}),(F_n,\phi_n)=(G,\psi)\]
in $\C$ such that for each $0\leq i\leq n-1$ there is a morphism $(F_i,\phi_i)\to(F_{i+1},\phi_{i+1})$ or a morphism $(F_{i+1},\phi_{i+1})\to(F_{i},\phi_{i})$. Recall that we have chosen a universal covering $(T,\tau)$ of $F$. Set $\eta_0:=\tau$. By inductively applying Lemma \ref{lemrepell} we obtain morphisms $\eta_i:(T,\phi_{i-1}\circ \eta_{i-1})\to (F_{i},\phi_{i})~(1\leq i\leq n)$ such that $(T,\eta_i)$ is a universal covering of $F_i$ for any $1\leq i\leq n$. Since the diagram
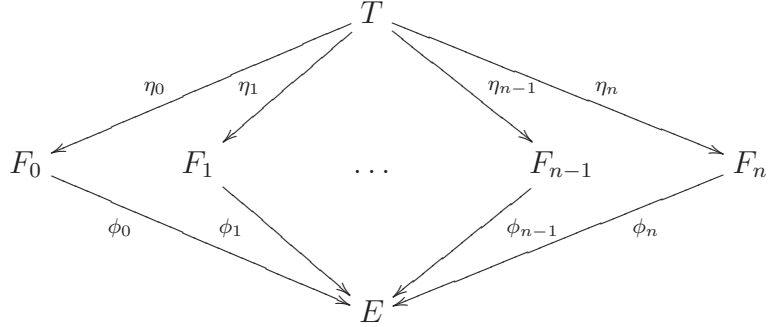
\[\xymatrix@C=45pt@R=40pt{&&T\ar[dll]_(0.6){\eta_0}\ar[dl]_(0.6){\eta_1}\ar[dr]^(0.6){\eta_{n-1}}\ar[drr]^(0.6){\eta_n}&&\\F_0\ar[drr]_(0.3){\phi_0}&F_1\ar[dr]_(0.3){\phi_1}&\dots&F_{n-1}\ar[dl]^(0.3){\phi_{n-1}}&F_n\ar[dll]^(0.3){\phi_n}\\&&E&&}\]
is commutative, we obtain that $\eta_n:(T_\C,\zeta_\C)=(T,\phi_0\circ\eta_0)\to (F_n,\phi_n)=(G,\psi)$ is a morphism in $\RG(\dot E)$. Thus $(T_\C,\zeta_\C)$ is a repelling object of $\C$. The second statement now follows from Proposition \ref{thmadm}.
\end{proof}

\begin{theorem}\label{thmm2}
Let $\C$ be a connected component of $\RG(\dot E)$. Then $(S_\C,\xi_\C)$ is an attracting object of $\C$, and consequently the objects of $\C$ are precisely the representation graphs $(G,\xi_\C\circ \alpha)$ where $(G,\alpha)$ is a covering of $S_\C$. 
\end{theorem}
\begin{proof}
The first statement of the theorem follows from Lemma \ref{lemcompare} and Theorem \ref{thmm1}. The second statement now follows from Propositions \ref{propsurj} and \ref{propsurj2}. 
\end{proof}

\begin{corollary}\label{corm2}
Let $\C$ be a connected component of $\RG(\dot E)$. Then up to isomorphism $(S_\C,\xi_\C)$ is the unique irreducible representation graph in $\C$. 
\end{corollary}
\begin{proof}
First we show that $(S_\C,\xi_\C)$ is irreducible. Suppose that $\xi_\C({}_{u}\!\Walk(S_\C))=\xi_\C({}_v\!\Walk(S_\C))$ where $u,v\in S_\C^0$. Recall that $(S_\C,\xi_\C)=((T_\C)_\sim,(\zeta_\C)_\sim)$. Hence $u=[x]_\sim$ and $v=[y]_\sim$ for some $x,y\in T_\C^0$. Let $\alpha: (T_\C,\zeta_\C)\to(S_\C,\xi_\C)$ be the morphism defined by $\alpha^0(w)=[w]_\sim$ and $\alpha^1(f)=[f]_\sim$ for any $w\in T_\C^0$ and $f\in T_\C^1$. It follows from Lemma \ref{propcruc} that 
\[\zeta_\C({}_{x}\!\Walk(T_\C))=\xi_\C({}_{u}\!\Walk(S_\C))=\xi_\C({}_v\!\Walk(S_\C))=\zeta_\C({}_{y}\!\Walk(T_\C)).\] 
Hence $x\sim y$ and thus $u=v$. Thus $(S_\C,\xi_\C)$ is irreducible.

Let now $(G,\psi)$ be an irreducible representation graph in $\C$. It follows from Proposition \ref{thmadm} and Theorem \ref{thmm2} that $(S_\C,\xi_\C)$ is isomorphic to a quotient of $(G,\psi)$. But since $(G,\psi)$ is irreducible, there is only one admissible equivalence relation on $G^0$, namely the equality relation $=$, and the corresponding quotient $(G_{=},\psi_{=})$ is isomorphic to $(G,\psi)$.
\end{proof}

\begin{example}\label{excatone}
Suppose that $E$ is the graph 
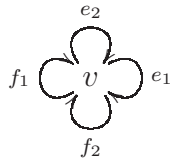
\[\xymatrix{ v\ar@(ur,dr)^{e_1}\ar@(ul,ur)^{e_2}\ar@(dl,ul)^{f_1}\ar@(dr,dl)^{f_2}}\]
and that $C=\{\{e_1,f_1\},\{e_2,f_2\}\}$ and $D=\{\{e_1,e_2\},\{f_1,f_2\}\}$. Consider the representation graphs $(F_1,\phi_1),\dots, (F_7,\phi_7)$ for $\dot E$ given below.

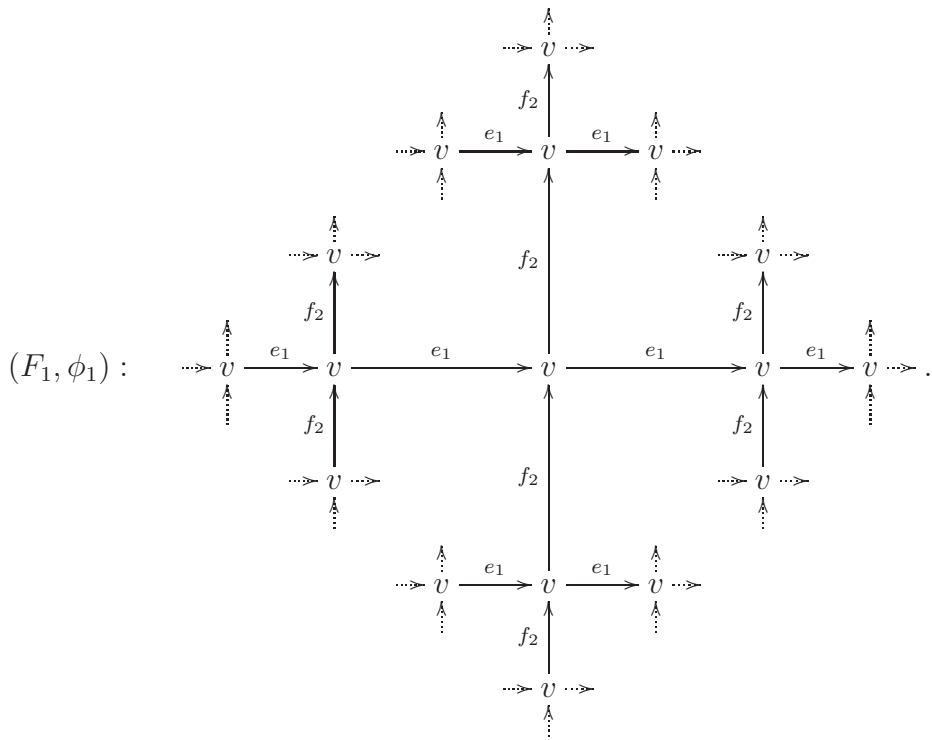
\begin{tabular}{l}
$\xymatrix@C=10pt@R=10pt{
&&&&&&&&&&&&&&&\\
&&&&&&&\ar@{..>}[r]&v\ar@{..>}[u]\ar@{..>}[r]&&&&&&&\\
&&&&&&&&&&&&&&&\\
&&&&&\ar@{..>}[r]&v\ar@{..>}[u]\ar^{e_1}[rr]&&v\ar[uu]^{f_2}\ar^{e_1}[rr]&&v\ar@{..>}[u]\ar@{..>}[r]&&&&&\\
&&&&&&\ar@{..>}[u]&&&&\ar@{..>}[u]&&&&&\\
&&&\ar@{..>}[r]&v\ar@{..>}[r]\ar@{..>}[u]&&&&&&&\ar@{..>}[r]&v\ar@{..>}[u]\ar@{..>}[r]&&&\\
&&&&&&&&&&&&&&&\\
(F_1,\phi_1):&\ar@{..>}[r]&v\ar@{..>}[u]\ar^{e_1}[rr]&&v\ar^{f_2}[uu]\ar^{e_1}[rrrr]&&&&v\ar^{f_2}[uuuu]\ar^{e_1}[rrrr]&&&&v\ar^{f_2}[uu]\ar^{e_1}[rr]&&v\ar@{..>}[u]\ar@{..>}[r]&.\\
&&\ar@{..>}[u]&&&&&&&&&&&&\ar@{..>}[u]&\\
&&&\ar@{..>}[r]&v\ar^{f_2}[uu]\ar@{..>}[r]&&&&&&&\ar@{..>}[r]&v\ar^{f_2}[uu]\ar@{..>}[r]&&&\\
&&&&\ar@{..>}[u]&&&&&&&&\ar@{..>}[u]&&&\\
&&&&&\ar@{..>}[r]&v\ar@{..>}[u]\ar^{e_1}[rr]&&v\ar^{f_2}[uuuu]\ar^{e_1}[rr]&&v\ar@{..>}[u]\ar@{..>}[r]&&&&&\\
&&&&&&\ar@{..>}[u]&&&&\ar@{..>}[u]&&&&&\\
&&&&&&&\ar@{..>}[r]&v\ar@{..>}[r]\ar^{f_2}[uu]&&&&&&&\\
&&&&&&&&\ar@{..>}[u]&&&&&&&
}$
\end{tabular}
\\
$~$
\vspace{0.2cm}
$~$
\\
\begin{tabular}{l}
$
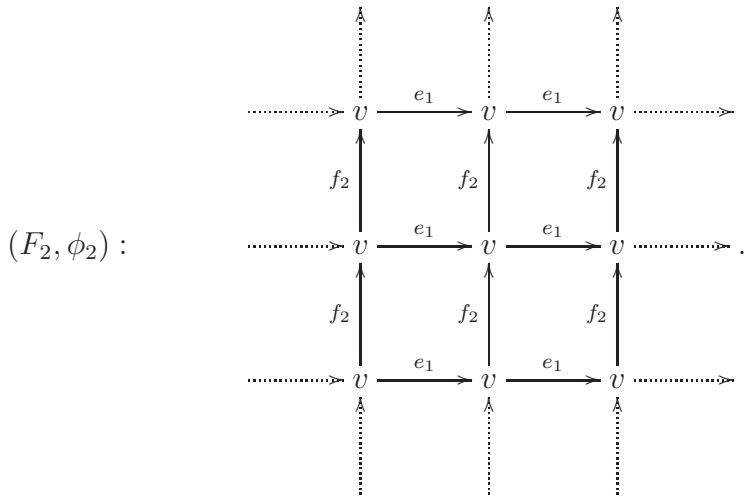
\xymatrix@C=35pt@R=35pt{
&&&&&\\
&\ar@{..>}[r]&v\ar^{e_1}[r]\ar@{..>}[u]&v\ar^{e_1}[r]\ar@{..>}[u]&v\ar@{..>}[r]\ar@{..>}[u]&\\
\hspace{1.1cm}(F_2,\phi_2):&\ar@{..>}[r]&v\ar^{e_1}[r]\ar^{f_2}[u]&v\ar^{e_1}[r]\ar^{f_2}[u]&v\ar^{f_2}[u]\ar@{..>}[r]&.\\
&\ar@{..>}[r]&v\ar^{e_1}[r]\ar^{f_2}[u]&v\ar^{e_1}[r]\ar^{f_2}[u]&v\ar^{f_2}[u]\ar@{..>}[r]&\\
&&\ar@{..>}[u]&\ar@{..>}[u]&\ar@{..>}[u]&
}$
\end{tabular}
\\
$~$
\vspace{0.5cm}
$~$
\\
\begin{tabular}{ll}
$
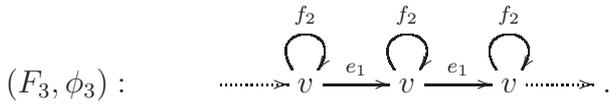
\xymatrix@C=25pt@R=25pt{
(F_3,\phi_3):&\ar@{..>}[r]&v\ar@(ul,ur)^{f_2}\ar^{e_1}[r]&v\ar@(ul,ur)^{f_2}\ar^{e_1}[r]&v\ar@(ul,ur)^{f_2}\ar@{..>}[r]&.}$
&$
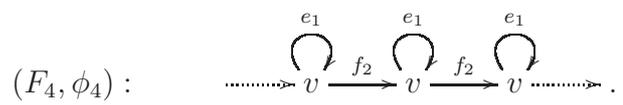
\xymatrix@C=25pt@R=25pt{
(F_4,\phi_4):&\ar@{..>}[r]&v\ar@(ul,ur)^{e_1}\ar^{f_2}[r]&v\ar@(ul,ur)^{e_1}\ar^{f_2}[r]&v\ar@(ul,ur)^{e_1}\ar@{..>}[r]&.
}$\\
\\
$
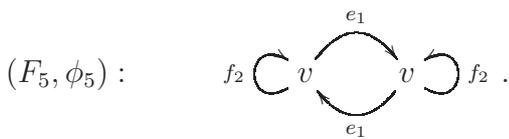
\xymatrix@C=25pt@R=25pt{
(F_5,\phi_5):&&v\ar@(dl,ul)^{f_2}\ar@/^1.3pc/^{e_1}[r]&v\ar@(dr,ur)_{f_2}\ar@/^1.3pc/^{e_1}[l]&.}$
&$
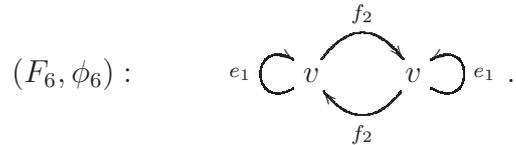
\xymatrix@C=25pt@R=25pt{
(F_6,\phi_6):&&v\ar@(dl,ul)^{e_1}\ar@/^1.3pc/^{f_2}[r]&v\ar@(dr,ur)_{e_1}\ar@/^1.3pc/^{f_2}[l]&.}$\\
\\
$
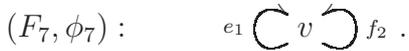
\xymatrix@C=25pt@R=25pt{
(F_7,\phi_7):&&v \ar@(dl,ul)^{e_1} \ar@(dr,ur)_{f_2}&.}$&
\end{tabular}\\
\\
\\
There are morphisms
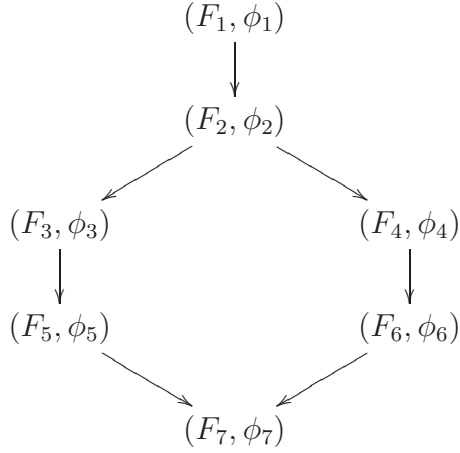
\[\xymatrix@C=20pt@R=20pt{
&(F_1,\phi_1)\ar[d]&\\
&(F_2,\phi_2)\ar[dl]\ar[dr]&\\
(F_3,\phi_3)\ar[d]&&(F_4,\phi_4)\ar[d]\\
(F_5,\phi_5)\ar[dr]&&(F_6,\phi_6)\ar[dl]\\
&(F_7,\phi_7)&
}\]
in $\RG(\dot E)$. Hence the representation graphs $(F_i,\phi_i)~(1\leq i\leq 7)$ lie in the same connected component $\C$ of $\RG(\dot E)$. We leave it to the reader to check that $(F_1,\phi_1)\cong (T_\C,\zeta_\C)$ and $(F_7,\phi_7)\cong (S_\C,\xi_\C)$.
\end{example}

\section{Modules for Leavitt path algebras of bi-separated graphs via representation graphs}\label{weightedrbgh}
Throughout this section $\dot E=(E,C,D)$ denotes a fixed bi-separated graph.

\subsection{The functor $V$}\label{vrepsect31}
For any object $(F,\phi)$ in $\RG(\dot E)$ let $V(F,\phi)$ be the $K$-vector space with basis $F^0$. 
For any $v\in E^0$ and $e\in E^1$ define endomorphisms $\sigma_v,\sigma_{e},\sigma_{e^*}\in \End_K(V(F,\phi))$ by
\begin{align*}
\sigma_{v}(w)     &=\begin{cases}w,\quad \quad                     &\text{if }\phi(w)=v,\\
0,     & \text{otherwise,} \end{cases}\\
\sigma_{e}(w)     &=\begin{cases}r(e),\quad                       &\text{if }\exists f\in F^1\text{ such that }s(f)=w\text{ and }\phi(f)=e,\\
0,     & \text{otherwise,} \end{cases}\\
\sigma_{e^*}(w)  &=\begin{cases}s(f),\quad                        &\text{if }\exists f\in F^1\text{ such that }r(f)=w\text{ and }\phi(f)=e,\\
0,    & \text{otherwise}, \end{cases}
\end{align*}
where $w\in F^0$. Note that $\sigma_{e}$ and $\sigma_{e^*}$ are well-defined since for any $w\in F^0$ the maps $\phi|_{s^{-1}(w)}$ and $\phi|_{r^{-1}(w)}$ are injective by Lemma \ref{lemwell}. It follows from the universal property of $L(\dot E)$ that there is a unique $K$-algebra homomorphism $L(\dot E)\rightarrow \End_K(V(F,\phi))^\op$ such that $\pi(v)=\sigma_{v}$, $\pi(e)=\sigma_e$ and $\pi(e^*)=\sigma_{e^*}$ for any $v\in E^0$ and $e\in E^1$. 
 
The vector space $V(F,\phi)$ becomes a right $L(\dot E)$-module by defining $x.a:= \pi(a)(x)$ for any $a\in L(\dot E)$ and $x\in V(F,\phi)$. We call $V(F,\phi)$ the {\it $L(\dot E)$-module defined by $(F,\phi)$}. A morphism $\alpha:(F,\phi)\rightarrow(G,\psi)$ in $\RG(\dot E)$ induces a surjective $L(\dot E)$-module homomorphism $V(\alpha):V(F,\phi)\rightarrow V(G,\psi)$ such that $V(\alpha)(u)=\alpha(u)$ for any $u\in F^0$. We obtain a functor
\[V:\RG(\dot E)\to \Modd(L(\dot E))\]
where $\Modd(L(\dot E))$ denotes the category of right $L(\dot E)$-modules.

Recall that any element of $L(\dot E)$ can be written as a $K$-linear combination of walks in $E$. The following lemma describes the action of $\Walk(E)$ on $V(F,\phi)$. Note that by Lemma \ref{lemwell}, for any $p\in \Walk(E)$ and $w\in F^0$ there is at most one walk $q\in \Walk (F)$ such that $s(q)=w$ and $\phi(q)=p$.

\begin{lemma}\label{lemaction}
Let $(F,\phi)$ be an object in $\RG(\dot E)$. If $p\in \Walk(E)$ and $w\in F^0$, then
\[w.p=\begin{cases}r(q),\quad\quad            &\text{if }\exists q\in \Walk(F) \text{ such that }s(q)=w\text{ and }\phi(q)=p,\\
0,                                           & \text{otherwise.} \end{cases}\]
\end{lemma}


The corollory below follows from Lemmas \ref{lemwell} and \ref{lemaction}.
\begin{corollary}\label{corwell}
If $p\in \Walk(E)$ and $w\neq w'\in F^0$, then either $w.p=w'.p=0$ or $w.p\neq w'.p$. 
\end{corollary}

\subsection{Fullness of the functor $V$}

Let $(F,\phi)$ be an object of $\RG(\dot E)$ and $(F_\approx,\phi_\approx)$ a quotient of $(F,\phi)$. By Proposition \ref{thmadm} there is a morphism $\phi:(F,\phi)\to (F_\approx,\phi_\approx)$ and hence a surjective morphism $V(\phi):V(F,\phi)\rightarrow V(F_\approx,\phi_\approx)$. By the lemma below, which is easy to check, there is also a morphism $V(F_\approx,\phi_\approx)\to V(F,\phi)$.

\begin{lemma}\label{lemopp}
Let $(F,\phi)$ be an object of $\RG(\dot E)$ and $(F_\approx,\phi_\approx)$ a quotient of $(F,\phi)$. Then there is a morphism $V(F_\approx,\phi_\approx)\to V(F,\phi)$ mapping $[u]\mapsto \sum_{v\approx u}v$.
\end{lemma}

In general $V$ is not full as the following example shows. Let $(F,\phi)$ and $(G,\psi)$ be finite representation graphs for $\dot E$ such that $(G,\psi)$ is a proper quotient of $(F,\phi)$ (for example $(F,\phi)=(F_6,\phi_6)$ and $(G,\psi)=(F_7,\phi_7)$ in Example \ref{excatone}). Then $F$ has strictly more vertices than $G$. By Lemma \ref{lemopp} there is a morphism $V(G,\psi)\to V(F,\phi)$. This morphism cannot be induced by a morphism $(G,\psi)\to (F,\phi)$. Otherwise $(F,\phi)$ would be isomorphic to a quotient of $(G,\psi)$, but this is not possible since $F$ has strictly more vertices than $G$. 

\begin{question}\label{Q1}
Do non-isomorphic representation graphs in $\RG(\dot E)$ define non-isomorphic modules for $L(\dot E)$?
\end{question}

The author does not know the answer to Question \ref{Q1}. But we will show that if $V(F,\phi)$ and $V( G,\psi)$ are isomorphic, then $(F,\phi)$ and $( G,\psi)$ lie in the same connected component of $\RG(\dot E)$.

\begin{lemma}\label{lemmodule}
Let $(F,\phi)$ and $( G,\psi)$ be objects in $\RG(\dot E)$ and $\theta:V(F,\phi)\rightarrow V( G,\psi)$ a morphism in $\Modd(L(\dot E))$. Suppose that $\theta(u)=\sum_{i=1}^nk_iv_i$ for some $u\in F^0$, $n\geq 1$, $k_1,\dots,k_n\in K^\times$ and pairwise distinct vertices $v_1,\dots,v_n\in G^{0}$. Then $\phi({}_u\!\Walk(F))=\psi({}_{v_i}\!\Walk( G))$ for any $1\leq i\leq n$.
\end{lemma}
\begin{proof} 
Let $p\in \Walk(\dot E)$ such that $p\not\in\phi({_u}\!\Walk(F))$. Then
\[0=\theta(0)=\theta(u.p)=\theta(u).p=(\sum_{i=1}^nk_iv_i).p=\sum_{i=1}^nk_i(v_i.p)\]
by Lemma \ref{lemaction}. It follows from Corollary \ref{corwell} that $v_i.p=0$ for any $1\leq i\leq n$, whence $p\not\in \psi({}_{v_i}\!\Walk( G))$ for any $1\leq i\leq n$. Hence we have shown that $\phi({}_u\!\Walk(F))\supseteq\psi({}_{v_i}\!\Walk( G))$ for any $1\leq i\leq n$. It follows from Lemma \ref{lempath} that $\phi({}_u\!\Walk(F))=\psi({}_{v_i}\!\Walk( G))$ for any $1\leq i\leq n$.
\end{proof}

 We leave the proof of the next lemma to the reader.

\begin{lemma}\label{lemwithoutnumber}
Let $(F,\phi)$ and $( G,\psi)$ be objects in $\RG(\dot E)$. Let $u\in F^{0}$ and $v\in G^{0}$ such that $\phi({}_u\!\Walk(F))=\psi({}_v\!\Walk( G))$. Then $\phi({}_{u.p}\!\Walk(F))=\psi({}_{v.p}\!\Walk( G))$ for any $p\in \phi({}_u\!\Walk(F))=\psi({}_v\!\Walk( G))$.
\end{lemma}

\begin{lemma}\label{prop41.5}
Let $(F,\phi)$ and $( G,\psi)$ be objects in $\RG(\dot E)$. Suppose that there exist  $u\in F^{0}$ and $v\in G^{0}$ such that $\phi({}_u\!\Walk(F))=\psi({}_v\!\Walk( G))$. Then $(F,\phi)$ and $( G,\psi)$ lie in the same connected component of $\RG(\dot E)$.
\end{lemma}
\begin{proof}
Since $F$ is connected, we can choose for any $u'\in F^{0}$ a $p_{u'}\in {}_{u}\!\Walk_{u'}(F)$. Define a graph homomorphism $\alpha:F_\sim\to G_\sim$ by 
\begin{align*}
&\alpha^0([u'])=[v.\phi(p_{u'})] \text{ for any }u'\in F^{0}\text{ and }\\
&\alpha^1([f])=[g]\text{ for any }f\in F^{1},\text{ where }s(g)=v.\phi(p_{s(f)})\text{ and }\psi(g)=\phi(f).
\end{align*}
It follows from Lemma \ref{lemwithoutnumber} that $\alpha$ is well-defined. One checks routinely that $\alpha:(F_\sim,\phi_\sim)\to( G_\sim,\psi_\sim)$ is a morphism in $\RG(\dot E)$. Thus, in view of Proposition \ref{thmadm}, $(F,\phi)$ and $( G,\psi)$ lie in the same connected component of $\RG(\dot E)$.
\end{proof}

The theorem below follows directly from Lemmas \ref{lemmodule} and \ref{prop41.5}.
\begin{theorem}\label{prop42}
Let $(F,\phi)$ and $( G,\psi)$ be objects of $\RG(\dot E)$. If there is a nonzero $L(\dot E)$-module homomorphism $\theta:V(F,\phi)\rightarrow V( G,\psi)$, then $(F,\phi)$ and $( G,\psi)$ lie in the same connected component of $\RG(\dot E)$.
\end{theorem}

\begin{corollary}\label{prop43}
Let $(F,\phi)$ and $( G,\psi)$ be irreducible representation graphs for $\dot E$. Then $(F,\phi)\cong( G,\psi)$ if and only if $V(F,\phi)\cong V( G,\psi)$.
\end{corollary}
\begin{proof}
Clearly $(F,\phi)\cong( G,\psi)$ implies $V(F,\phi)\cong V( G,\psi)$ since $V$ is a functor. Suppose now that $V(F,\phi)\cong V( G,\psi)$. Then, by Theorem \ref{prop42}, $(F,\phi)$ and $( G,\psi)$ lie in the same connected component $\C$ of $\RG(\dot E)$. It follows from Corollary \ref{corm2} that $(F,\phi)\cong (S_\C,\xi_\C)\cong ( G,\psi)$.
\end{proof}

\subsection{Simplicity of the modules $V(F,\phi)$}
In this subsection we show that the $L(\dot E)$-module $V(F,\phi)$ is simple if and only if $(F,\phi)$ is irreducible.

\begin{lemma}[{\cite[Lemma B.6]{HPS}}]\label{lembasis}
Let $W$ be a $K$-vector space and $B$ a linearly independent subset of  $W$. Let $k_i\in K$ and $u_i,v_i\in B$, where $1\leq i \leq n$. Then $\sum_{i=1}^nk_i(u_i-v_i)\not\in B$.
\end{lemma}


\begin{theorem}\label{thmirr} Let $(F,\phi)$ be an object in $\RG(\dot E)$
Then the following are equivalent.
\begin{enumerate}[\upshape(i)]
\item $V(F,\phi)$ is simple.
\medskip
\item For any $x\in V(F,\phi)\setminus\{0\}$ there is an $a\in L(\dot E)$ such that $x.a \in F^{0}$.
\medskip
\item For any $x\in V(F,\phi)\setminus\{0\}$ there is a $k\in K$ and a $p\in \Walk(E), $ such that $x.kp\in F^{0}$.
\medskip
\item $(F,\phi)$ is irreducible.
\end{enumerate}
\end{theorem}
\begin{proof} 
(i) $\Rightarrow$ (iv). Assume that there are $u\neq v\in F^{0}$ such that $\phi({}_u\!\Walk(F))= \phi({}_v\!\Walk(F))$. Consider the submodule $(u-v).L(\dot E)\subseteq V(F,\phi)$. Since $V(F,\phi)$ is simple by assumption, we have $(u-v).L(\dot E)=V(F,\phi)$. Hence there is an $a\in L(\dot E)$ such that $(u-v).a=v$. Clearly there is an $n\geq 1$, $k_1,\dots,k_n\in K^\times$ and pairwise distinct $p_1,\dots,p_n\in \Walk(E)$ such that $a=\sum_{i=1}^nk_ip_i$. We may assume that $(u-v).p_i\neq 0$ for any $1\leq i\leq n$. It follows that $p_i\in \phi({}_u\!\Walk(F))= \phi({}_v\!\Walk(F))$ for any $i$ and moreover, that $(u-v).p_i=u_i-v_i$ for some distinct $u_i,v_i\in F^{0}$. Hence
\[v=(u-v).a=(u-v).(\sum_{i=1}^nk_ip_i)=\sum_{i=1}^nk_i(u_i-v_i)\]
which contradicts Lemma \ref{lembasis}.

\medskip 
(iv) $\Rightarrow$ (iii). Let $x\in V(F,\phi)\setminus \{0\}$. Then there is an $n\geq 1$, $k_1,\dots,k_n\in K^\times$ and pairwise distinct $v_1,\dots, v_n\in F^{0}$ such that $x=\sum_{i=1}^nk_iv_i$. If $n=1$, then $x.k_1^{-1}\phi(v_1) = v_1$. Suppose now that $n>1$. By assumption, we can choose a $p_1\in \phi({}_{v_1}\!\Walk(F))$ such that $p_1\not\in\phi({}_{v_2}\!\Walk(F))$. Clearly $x.p_1\neq 0$ is a linear combination of at most $n-1$ vertices from $F^{0}$. Proceeding this way, we obtain walks $p_1,\dots,p_m$ such that $x.p_1\dots p_m=kv$ for some $k\in K^\times$ and $v\in F^{0}$. Hence $x.k^{-1}p_1\dots p_m=v$.

\medskip 

(iii) $\Rightarrow$ (ii). This implication is trivial.

\medskip 

(ii) $\Rightarrow$ (i). Let $U\subseteq V(F,\phi)$ be a nonzero $L(\dot E)$-submodule and $x\in U\setminus\{0\}$. By assumption, there is an $a\in L(\dot E)$ and a $v\in F^{0}$ such that $v=x.a\in U$. Let now $v'$ be an arbitrary vertex in $F^{0}$. Since $F$ is connected, there is a $p\in {}_{v}\!\Walk_{v'}(F)$. It follows that $v'=v.\phi(p)\in U$. Hence $U$ contains $F^{0}$ and thus $U=V(F,\phi)$.
\end{proof}

\begin{remark}
If the bi-separation on $E$ is the Cuntz-Krieger bi-separation and hence $L(\dot{E})$ is isomorphic to the usual Leavitt path algebra $L(E)$ (cf. \cite{mohan-suhas}), then the simple modules $V(F,\phi)$ where $(F,\phi)$ is an irreducible representation graph for $\dot E$ are precisely the Chen simple modules (cf. \cite[\S 4]{HPS}).
\end{remark}

\section{Indecomposability of the modules $V(T_\C,\zeta_{\C})$}
In this section $\dot E$ denotes a fixed bi-separated graph and $\C$ a connected component of $\RG(\dot E)$. Recall that for a ring $R$, an $R$-module is called \emph{indecomposable} if it is nonzero and cannot be written as a direct sum of two nonzero submodules. It is easy to see that an $R$-module $M$ is indecomposable if and only if $\End_R(M)$ has no {\it nontrivial} idempotents, i.e. idempotents distinct from $0$ and $1$. We will show that the $L(\dot E)$-module $V(T_\C,\zeta_{\C})$ is indecomposable. In order to do so, we first prove some technical lemmas. Then we use these lemmas to prove Propositions \ref{propextend}, \ref{propfree} and \ref{propKG}. Finally we deduce the main result of this section, Theorem \ref{thmdecomp}, from these propositions.

In the following we may write $( T,\zeta)$ instead of $( T_{\C},\zeta_{\C})$ and $( S,\xi)$ instead of $( S_{\C},\xi_{\C})$. By Theorem \ref{thmm2} there is a morphism $\eta:(T,\zeta)\to( S,\xi)$. Hence the diagram
\[\xymatrix{
T\ar[r]^\eta\ar@/^2pc/[rr]^\zeta&S\ar[r]^\xi&E
}\]
commutes. By Proposition \ref{propsurj}, $(T,\eta)$ is a covering of $S$. We fix a $y\in  S^0$ and denote the linear subspace of $V( T,\zeta)$ with basis $\eta^{-1}(y)$ by $W$. Moreover, we denote the subalgebra of $L(\dot E)$ consisting of all $K$-linear combination of elements of $\xi({}_y\!\Walk_y( S))$ by $A$ (note that $\xi({}_y\!\Walk_y( S))\subseteq \Walk(E)$). One checks easily that the action of $L(\dot E)$ on $V( T,\zeta)$ induces an action of $A$ on $W$, making $W$ a right $A$-module. Set $\bar A:=A/\Ann(W)$ and let $~\bar~:A\to \bar A,~a\mapsto \bar a$ be the canonical algebra homomorphism. The action of $A$ on $W$ induces an action of $\bar A$ on $W$ making $W$ a right $\bar A$-module. 

If $p\in \xi({}_y\!\Walk_y( S))$, then there is a unique $\hat p\in {}_y\!\Walk_y( S)$ such that $\xi(\hat p)=p$ (the uniqueness follows from Lemma \ref{lemwell}). Since $( T,\eta)$ is a covering of $ S$, there is for any $x\in \eta^{-1}(y)$ a unique $\tilde p_x\in {}_x\!\Walk( T)$ such that $\eta(\tilde p_x)=\hat p$. We define a semigroup homomorphism $f:\xi({}_y\!\Walk_y( S))\to \pi( S,y)$ by $f(p)=\underline{\hat p}$.

\begin{lemma}\label{lemcous}
Let $x\in  T^{0}$ and $p,q\in {}_x\!\Walk (T)$. Then $r(p)=r(q)$ if and only if $\underline{\eta(p)}=\underline{\eta(q)}$.
\end{lemma}
\begin{proof}
($\Rightarrow$) Suppose that $r(p)=r(q)$. Then $pq^*\in {}_x\!\Walk_x( T)$. Since $T$ is a tree it follows that $\underline{pq^*}=x$. Hence $(\underline{\eta(p)})(\underline{\eta(q)^*})=\underline{\eta(p)\eta(q)^*}=\eta(\underline{pq^*})=\eta(x)$ by Lemma \ref{lemnospurs} and thus $\underline{\eta(p)}=\underline{\eta(q)}$.\\
($\Leftarrow$). Suppose that $\underline{\eta(p)}=\underline{\eta(q)}$. Then $\eta(\underline{p})=\eta(\underline{q})$ by Lemma \ref{lemnospurs}. Since $( T,\eta)$ is a covering of $ S$, it follows that $\underline{p}=\underline{q}$. Thus $r(p)=r(q)$.
\end{proof}

The lemma below follow from Lemma \ref{lemcous}.

\begin{lemma}\label{lemzedp}
Let $p,q\in \xi({}_y\!\Walk_y( S))$. Then the following are equivalent.
\begin{enumerate}[(i)]
\item $f(p)=f(q)$.
\item $r(\tilde p_x)=r(\tilde q_x)$ for some $x\in\eta^{-1}(y)$.
\item $r(\tilde p_x)=r(\tilde q_x)$ for any $x\in\eta^{-1}(y)$.
\end{enumerate}
\end{lemma}

We set $\Ann(x):=\{a\in A\mid x.a=0\}$ for any $x\in\eta^{-1}(y)$, and $\Ann(W):=\{a\in A\mid W.a=0\}$.

\begin{lemma}\label{lemann}
Let $x\in\eta^{-1}(y)$. Then 
\[\Ann(x)=\{\sum_{p\in\xi({}_y\!\Walk_y( S))}k_pp\in A\mid \sum_{\substack{p\in \xi({}_y\!\Walk_y( S)),\\r(\tilde p_x)=x'}}k_p=0 \quad \text{for any } x'\in \eta^{-1}(y)\}\]
\end{lemma}
\begin{proof}
Let $a=\sum_{p\in\xi({}_y\!\Walk_y( S))}k_pp\in A$. Then, in view of Lemma \ref{lemaction},
\begin{align*}
&\hspace{2.35cm}a\in\Ann(x)&\\
&\hspace{1.45cm}\Leftrightarrow~x.\sum_{p\in\xi({}_y\!\Walk_y( S))}k_pp=0&
\end{align*}
\begin{align*}
\Leftrightarrow~&\sum_{x'\in\eta^{-1}(y)}(\sum_{\substack{p\in \xi({}_y\!\Walk_y( S)),\\r(\tilde p_x)=x'}}k_p)x'=0\\
\Leftrightarrow~&\sum_{\substack{p\in \xi({}_y\!\Walk_y( S)),\\r(\tilde p_x)=x'}}k_p=0 \quad\forall x'\in \eta^{-1}(y).
\end{align*}
\end{proof}

\begin{corollary}\label{corann}
$\Ann(W)=\Ann(x)$ for any $x\in\eta^{-1}(y)$.
\end{corollary}
\begin{proof}
Since $\Ann(W)=\bigcap_{x\in\eta^{-1}(y)}\Ann(x)$, it suffices to show that $\Ann(x_1)=\Ann(x_2)$ for any $x_1,x_2\in \eta^{-1}(y)$. So let $x_1,x_2\in \eta^{-1}(y)$. For any $x\in \eta^{-1}(y)$ set $Y_{x}:= \{p\in \xi({}_y\!\Walk_y( S))\mid r(\tilde p_{x_1})=x\}$ and $Z_{x}:= \{p\in \xi({}_y\!\Walk_y( S))\mid r(\tilde p_{x_2})=x\}$. Clearly $\xi({}_y\!\Walk_y( S))=\bigsqcup_{x\in\eta^{-1}(y)} Y_x=\bigsqcup_{x\in\eta^{-1}(y)} Z_x$. It follows from Lemma \ref{lemzedp} that there is a permutation $\pi$ on the set $\eta^{-1}(y)$ such that $Y_x=Z_{\pi(x)}$ for any $x\in\eta^{-1}(y)$. Let now $a=\sum_{p\in\xi({}_y\!\Walk_y( S))}k_pp\in A$. Then, by Lemma \ref{lemann}, 
\begin{align*}
&a\in\Ann(x_1)\\
\Leftrightarrow~&\sum_{p\in Y_x} k_p=0~\forall x\in\eta^{-1}(y)\\
\Leftrightarrow~&\sum_{p\in Z_x} k_p=0~\forall x\in\eta^{-1}(y)\\
\Leftrightarrow~&a\in\Ann(x_2).
\end{align*}
\end{proof}

Recall that if $x,x'\in T^{0}$, then $x\sim x'\Leftrightarrow \zeta({}_x\!\Walk( T))=\zeta({}_{x'}\!\Walk( T))$.  

\begin{lemma}\label{lemalphsim}
$\eta^{-1}(y)$ is a $\sim$-equivalence class.
\end{lemma}
\begin{proof}
Choose an $x\in \eta^{-1}(y)$. Let $x'\in  T^{0}$ such that $x'\sim x$. It follows from Lemma \ref{lempath} that $\xi({}_{\eta(x)}\!\Walk( S))=\zeta({}_{x}\!\Walk( T))=\zeta({}_{x'}\!\Walk( T))=\xi({}_{\eta(x')}\!\Walk( S))$. Hence $\eta(x)\sim\eta(x')$. It follows that $\eta(x')=\eta(x)=y$ since $ S$ is irreducible.\\
Let now $x'\in \eta^{-1}(y)$. Then $\eta({}_x\!\Walk( S))=\eta({}_{x'}\!\Walk( S))$ since $( T,\eta)$ is a covering of $ S$. Hence $\zeta({}_x\!\Walk( S))=\zeta({}_{x'}\!\Walk( S))$, i.e. $x\sim x'$. We have shown that $\eta^{-1}(y)=[x]_{\sim}$.
\end{proof}

\begin{proposition}\label{propextend}
Any nontrivial idempotent endomorphism in $\End_{L(\dot E)}(V( T,\zeta))$ restricts to a nontrivial idempotent endomorphism in $\End_{\bar A}(W)$. 
\end{proposition}
\begin{proof}
Let $\epsilon\in \End_{L(\dot E)}(V_{( T,\zeta)})$ be a nontrivial idempotent endomorphism. It follows from Lemmas \ref{lemmodule} and \ref{lemalphsim} that $\epsilon(W)\subseteq W$. Hence $\epsilon|_W\in \End_{\bar A}(W)$. Clearly $\epsilon|_W$ is an idempotent. It remains to show that $\epsilon|_W$ is nontrivial. Since $\epsilon$ is nontrivial, there are $v,w\in  T^{0}$ such that $\epsilon(v)\neq 0$ and $\epsilon(w)\neq w$. Let $x\in\eta^{-1}(y)$ and choose a $p\in \zeta({}_x\!\Walk_v( T))$ and a $q\in \zeta({}_x\!\Walk_w( T))$. Then $\epsilon(v)=\epsilon( x.p)=\epsilon(x).p$ and $\epsilon(w)=\epsilon( x.q)=\epsilon(x).q$. It follows that $\epsilon(x)\neq 0, x$. Thus $\epsilon|_W$ is nontrivial.
\end{proof}

\begin{proposition}\label{propfree}
$W$ is free of rank $1$ as an $\bar A$-module.
\end{proposition}
\begin{proof}
Choose an $x\in\eta^{-1}(y)$. If $x'\in\eta^{-1}(y)$, then there is a $p\in{}_{x}\!\Walk_{x'}( T)$ since $ T$ is connected. Clearly $\zeta(p)\in \xi({}_y\!\Walk_y( S))$ and moreover $x.\overline{\zeta(p)}=x'$. Hence $x$ generates the $\bar A$-module $W$. On the other hand, if $x.\bar a=0$ for some $a\in A$, then $a\in \Ann(x)$ (since  $x.\bar a=x.a$) and hence $\bar a=0$ by Corollary \ref{corann}. Thus $\{x\}$ is a basis for the $\bar A$-module $W$.
\end{proof}

Recall that the semigroup homomorphism $f:\xi({}_y\!\Walk_y( S))\to \pi( S,y)$ was defined by $f(p)=\underline{\hat p}$. 

\begin{proposition}\label{propKG}
The algebra $\bar A$ is isomorphic to the group algebra $KG$ where $G=\pi( S,y)$.
\end{proposition}
\begin{proof}
Define the map
\begin{align*}
F:\bar A&\to KG,\\
\overline{\sum_{p\in\xi({}_y\!\Walk_y( S))}k_pp}&\mapsto \sum_{p\in\xi({}_y\!\Walk_y( S))}k_pf(p).
\end{align*}
First we show that $F$ is well-defined. Choose an $x\in\eta^{-1}(y)$. Suppose that $\overline{\sum_{p\in\xi({}_y\!\Walk_y( S))}k_pp}=\overline{\sum_{p\in\xi({}_y\!\Walk_y( S))}l_pp}$. Then $\sum_{p\in\xi({}_y\!\Walk_y( S))}(k_p-l_p)p\in \Ann(W)\subseteq \Ann(x)$. It follows from Lemma \ref{lemann} that 
\begin{equation}
\sum_{\substack{p\in \xi({}_y\!\Walk_y( S)),\\r(\tilde p_x)=x'}}(k_p-l_p)=0 \quad \text{for any } x'\in \eta^{-1}(y).
\end{equation}
We have to show that $\sum_{p\in\xi({}_y\!\Walk_y( S))}k_pf(p)=\sum_{p\in\xi({}_y\!\Walk_y( S))}l_pf(p)$, i.e.
\begin{equation}
\sum_{\substack{p\in \xi({}_y\!\Walk_y( S)),\\f(p)=g}}(k_p-l_p)=0 \quad \text{for any } g\in G.\quad\quad\hspace{0.05cm}{}
\end{equation}
But it follows from (1) and Lemma \ref{lemzedp} that (2) holds. Hence $F$ is well-defined. We leave it to the reader to check that $F$ is an algebra homomorphism. It remains to show that $F$ is bijective. Suppose that $F( \overline{\sum_{p\in\xi({}_y\!\Walk_y( S))}k_pp})=F(\overline{\sum_{p\in\xi({}_y\!\Walk_y( S))}l_pp})$. Then (2) holds. It follows from (2) and Lemma \ref{lemzedp} that (1) holds. Hence $\overline{\sum_{p\in\xi({}_y\!\Walk_y( S))}k_pp}=\overline{\sum_{p\in\xi({}_y\!\Walk_y( S))}l_pp}$ and therefore $F$ is injective. The surjectivity of $F$ follows from the surjectivity of $f$.
\end{proof}

We are now in position to prove the main result of this section.
\begin{theorem}\label{thmdecomp}
The $L(\dot E)$-module $V( T_{\C},\zeta_{\C})$ is indecomposable.
\end{theorem}
\begin{proof}
It follows from Propositions \ref{propextend}, \ref{propfree} and \ref{propKG} that
\begin{align*}
&V( T,\zeta)\text{ is indecomposable }\\
\Leftrightarrow~&\End_{L(\dot E)}(V( T,\zeta))\text{ has no nontrivial idempotents}\\
\Leftarrow~&\End_{\bar A}(W)\text{ has no nontrivial idempotents}\\
\Leftrightarrow~&\bar A\text{ has no nontrivial idempotents}\\
\Leftrightarrow~&KG\text{ has no nontrivial idempotents}
\end{align*}
where $G=\pi( S,y)$. But it is well-known that fundamental groups of graphs are free (see for example \cite[Lemma 4.10]{kp}). Hence the group ring $KG$ has no zero divisors by \cite[Theorem 12]{higman}. It follows that $KG$ has no nontrivial idempotents and thus $V( T_{\C},\zeta_{\C})$ is indecomposable.
\end{proof}

\section{Branching systems for bi-separated graphs} \label{branchhonda}
In this section $\dot E$ denotes a fixed bi-separated graph. We will introduce the notion of an $\dot E$-algebraic branching system and show that any $\dot E$-algebraic branching system induces a module for $L(\dot E)$. Then we will show that the categories $\RG(\dot E)$ of representation graphs for $\dot E$ and $\ABS(\dot E)$ of $\dot E$-algebraic branching systems are equivalent. In order to obtain this equivalence we no longer assume that representation graphs are connected!

\subsection{$\dot E$-algebraic branching systems}


\begin{definition}\label{defabs}
Let
\begin{itemize}
\item $\{Q_v\}_{v\in E^0}$ be a family of pairwise disjoint subsets of a nonempty set $Q$,
\item $\{S_{e}\}_{e\in E^1}$ a family of sets such that $Q_v=\bigsqcup_{e\in X}S_{e}$ for any $v\in E^0$ and $X\in C_v$,
\item $\{R_{e}\}_{e\in E^1}$ a family of sets such that $Q_v=\bigsqcup_{e\in Y}R_{e}$ for any $v\in E^0$ and $Y\in D_v$, and
\item $\{\rho_{e}\}_{e\in E^1}$ a family of bijections $\rho_{e}:S_{e} \to R_e$.
\end{itemize}
Then $(Q,\{Q_v\}_{v\in E^0},\{S_{e}\}_{e\in E^1},\{R_{e}\}_{e\in E^1},\{\rho_{e}\}_{e\in E^1})$
is called an {\it $\dot E$-algebraic branching system}. 
In this paper all $\dot E$-algebraic branching systems are assumed to be \textit{saturated}, i.e. $Q=\bigcup_{v\in E^0}Q_v$.
\end{definition}

We denote by $\ABS(\dot E)$ the category of $\dot E$-algebraic branching systems. A morphism $\alpha:\Es\to\Es'$ between two objects $\Es=(Q,\{Q_v\},\{S_{e}\},\{R_e\},\{\rho_{e}\})$ and $\Es'=(Q',\{Q'_v\},\{S'_{e}\},\{R'_e\},\{\rho'_{e}\})$ in $\ABS(\dot E)$ is a map $\alpha: Q\to Q'$ such that $\alpha(Q_v)\subseteq Q'_v$, $\alpha (S_{e})\subseteq S'_{e}$ and $\alpha(R_e)\subseteq R'_e$ for any $v\in E^0$ and $e\in E^1$, and $\alpha$ is compatible with the bijections inside $\Es$ and $\Es'$. 

Suppose that $\Es=(Q,\{Q_v\},\{S_{e}\},\{R_e\},\{\rho_{e}\})$ is an object in $\ABS(\dot E)$. Let $W=W(\Es)$ denote the $K$-vector space with basis $Q$. For any $v\in E^0$ and $e\in E^1$ define endomorphisms $\sigma_v,\sigma_e,\sigma_{e^*}\in \End_K(W)$ by 

\begin{align*}
\sigma_v(q)&=\begin{cases}
q,\quad\quad\hspace{0.12cm}&\text{ if }q\in Q_v,\\
0,&\text{ otherwise,}
\end{cases}\\
\sigma_{e}(q)&
=\begin{cases}
\rho_{e}(q),\hspace{0.24cm}&\text{ if }q\in S_{e},\\
0,&\text{ otherwise,}
\end{cases}\\
\sigma_{e^*}(q)&=\begin{cases}
\rho_{e}^{-1}(q),\hspace{0.0cm}&\text{ if }q\in R_{e},\\
0,&\text{ otherwise,}
\end{cases}
\end{align*}
for any $q\in Q$. It follows from the universal property of $L(\dot E)$ that there is a $K$-algebra homomorphism $\pi:L(E)\to \End_K(W)^{\op}$ such that $\pi(v)=\sigma_v$, $\pi(e)=\sigma_e$ and $\pi(e^*)=\sigma_{e^*}$ for any $v\in E^0$ and $e\in E^1$. 

If $\Es$ is an object in $\ABS(\dot E)$, then the vector space $W(\Es)$ becomes a right $L(\dot E)$-module by defining $w.q:=\pi(a)(w)$ for any $w\in W(\Es)$ and $a\in L(\dot E)$. We call this module the \textit{$L(\dot E)$-module defined by $\Es$}. If $\alpha:\Es\to\Es'$ is a morphism in $\ABS(\dot E)$, let $W(\alpha):W(\Es)\to W(\Es')$ be the module homomorphism such that $W(\alpha)(q)=\alpha(q)$ for any $q\in Q$. We obtain a functor $W:\ABS(\dot E)\to \Modd(L(\dot E))$.

\subsection{Branching systems vs. representation graphs}
To any object $\Es=(Q,\{Q_v\},\{S_{e}\},\{R_e\},$ $\{\rho_{e}\})$ in $\ABS (\dot E)$ we associate the object $\eta(\Es)=(F,\phi)$ in $\RG(\dot E)$ defined by
\begin{align*}
F^0&=Q,\\
F^1&=\{f_{q,X}\mid q\in Q_v\text{ and }X\in C_v \text{ for some }v\in E^0\},\\
s(f_{q,X})&=q,\\
r(f_{q,X})&=\rho_e(q)\text{ where }e\in X\text{ has the property that } q\in S_e,\\
\phi^0(q)&=v\text{ if }q\in Q_v,\\
\phi^1(f_{q,X})&=e\text{ where }e\in X\text{ has the property that } q\in S_e.
\end{align*}
To any morphism $\alpha:\Es=(Q,\{Q_v\},\{S_{e}\},\{R_e\},\{\rho_{e}\})\to\Es'=(Q',\{Q'_v\},\{S'_{e}\},\{R'_e\},\{\rho'_{e}\})$ in $\ABS (\dot E)$ we associate the morphism $ \eta(\alpha):\eta(\Es)\to \eta(\Es')$ in $\RG(\dot E)$ defined by $ \eta(\alpha)^0(q)=\alpha(q)$ for any $q\in Q$ and $ \eta(\alpha)^1(f_{q,X})=f'_{\alpha(q),X}$ for any $q\in Q_v$ and $X\in C_v$ where $v\in E^0$. In this way we obtain a functor $\eta:\ABS (\dot E)\to \RG(\dot E)$. 

Conversely, to any object $(F,\phi)$ in $\RG(\dot E)$ we associate an object $ \theta(F,\phi)=(Q,\{Q_v\},\{S_{e}\},\{R_e\},$ $\{\rho_{e}\})$ in $\ABS (\dot E)$ defined by
\begin{align*}
Q&=F^0,\\
Q_v&=\{w\in F^0\mid \phi(w)=v\}\text{ for any }v\in E^0,\\
S_{e}&=\{w\in F^0\mid \exists f\in s^{-1}(w):\phi(f)=e\}\text{ for any }e\in E^1,\\
R_{e}&=\{w\in F^0\mid \exists f\in r^{-1}(w):\phi(f)=e\}\text{ for any }e\in E^1,\\
\rho_{e}(w)&=r(f)\text{ for any }e\in E^1\text{ and }w\in S_e,\text{ where }f\in s^{-1}(w)\\\
&\quad\hspace{3.5cm}\text{ has the property that }\phi(f)=e. 
\end{align*}
To any morphism $\alpha:(F,\phi)\to (G,\psi)$ in $\RG(\dot E)$ we associate the morphism $ \theta(\alpha): \theta(F,\phi)\to \theta(G,\psi)$ in $\ABS(\dot E)$ defined by $ \theta(\alpha)(q)=\alpha(q)$ for any $q\in Q=F^0$. In this way we obtain a functor $ \theta:\RG(\dot E)\to \ABS (\dot E)$. 

We leave it to the reader to check that $ \theta\circ\eta=\id_{\ABS (\dot E)}$ and $ \eta\circ   \theta\cong\id_{\RG(\dot E)}$. Hence the categories $\RG(\dot E)$ and $\ABS (\dot E)$ are equivalent. Moreover, the diagrams
\[\xymatrix{\ABS (\dot E)\ar[dr]^W&\\&\Modd L(\dot E)& \text{and}\\\RG(\dot E)\ar[uu]^{\theta}\ar[ur]_V&\\}
\quad\quad
\xymatrix{\ABS (\dot E)\ar[dd]_{\eta}\ar[dr]^W&\\&\Modd L(\dot E)\\\RG(\dot E)\ar[ur]_V&\\}
\]
are commutative.

\end{document}